\documentclass[twoside]{article}

\usepackage{amssymb}
\usepackage{amsmath}
\usepackage{cite}
\usepackage{mathrsfs}

\makeatletter

\def\address#1{\gdef\@address{#1}}
\address{\bigskip\hfil\begin{tabular}{l@{}}
            School of Mathematics and Statistics F07\\  
            University of Sydney,
            Sydney N.S.W. 2006.\hfill\qquad
                        {\tt mathas@maths.usyd.edu.au}\\  
            Australia.\hfill\qquad
                    {\tt www.maths.usyd.edu.au:8000/u/mathas/}
          \end{tabular}}

\let\@classification\relax\let\@support\relax
\newif\ifClassificationOrSupport\ClassificationOrSupportfalse
\def\AMSClassification#1{\ClassificationOrSupporttrue
  \def\@classification{A.M.S. subject classification (1991): #1.\\ }}
\def\Support#1{\ClassificationOrSupporttrue\def\@support{#1}}

\def\@righthead{\ifnum\c@page=1\else%
       $\underline{\strut\hbox to\textwidth{\rlap{\bf\thepage}%
                   \hfil{\sc\@author}\hfil}}$\fi}
\def\@lefthead{\ifnum\c@page=1\else%
       $\underline{\strut\hbox to \textwidth{\hfil{\sc\@title}%
                    \hfil\llap{\bf\thepage}}}$\fi}

\def\ps@mine{
  \if@twoside\let\@oddhead=\@lefthead\let\@evenhead=\@righthead
  \else
    \def\@oddhead{\ifodd\c@page\@lefthead\else\@righthead\fi}
  \fi
  \def\@oddfoot{\hbox to \textwidth{%
                 \ifnum\c@page=1\hfil{\bf\thepage}\hfil\fi}}
}

\let\ps@plain\ps@mine
\let\atop\@@atop

\def\And{\text{\ and\ }}

\def\ForSome{\text{\ for some\ }}
\def\If{\text{\ if\ }}

\def\Whenever{\text{\ whenever\ }}
\def\Otherwise{\text{\ otherwise}}

\def\proof{\noindent{\bf Proof}\space\space}
\def\proofof#1{\noindent{\bf Proof of \ref{#1}\space}}

\newbox\squ  
\setbox\squ=\hbox{\vrule width.6pt
            \vbox{\hrule height.6pt width.4em\kern1ex\hrule height.6pt}%
            \vrule width.6pt}
\def\endproof{%
    \ifmmode\eqno\copy\squ\medskip\else{\unskip\nobreak\hfil%
    \penalty50\hskip2em\hbox{}\nobreak\hfil\copy\squ
    \parfillskip=0pt \finalhyphendemerits=0\penalty-100\medskip}\fi
  }

\def\Set[#1]#2|#3|{\Big\{\ #2\ \Big| \
            \vcenter{\hsize #1mm\centering#3}\Big\}}

{\catcode`\|=\active
  \gdef\set#1{\mathinner{\lbrace\,{\mathcode`\|"8000%
                                   \let|\midvert #1}\,\rbrace}}
}
\def\midvert{\egroup\mid\bgroup}

\def\Number#1{\refstepcounter{equation}\@last{(\theequation)}%
              \leqno(\theequation)\if*#1%
              \else\def\@currentlabel{{\rm(\theequation)}}\label{#1}%
              \fi}
\def\Dag{\ifmmode\leqno(\dag)\else$(\dag\)$\fi}
\def\DDag{\ifmmode\leqno(\ddag)\else$(\ddag\)$\fi}

\let\last\relax
\let\Last\relax
\def\@last#1{\global\let\last\Last\global\edef\Last{#1}}

\renewenvironment{enumerate}%
  {\ifnum\@enumdepth>3\@toodeep\else
     \advance\@enumdepth\@ne
     \edef\@enumctr{enum\romannumeral\the\@enumdepth}%
     \topsep\z@\parskip\z@
     \list{\csname label\@enumctr\endcsname}
       {\@nmbrlisttrue\let\@listctr\@enumctr
       \parsep\z@\itemsep\z@\topsep\z@
       \setcounter{\@enumctr}0
       \def\makelabel##1{\hss\llap{\rm ##1}}
     }
   \fi}
  {\endlist}

\def\thebibliography#1{%
  \section*{References\@mkboth{References}{References}}\list
  {[\arabic{enumi}]}{\settowidth\labelwidth{[#1]}\leftmargin\labelwidth
  \advance\leftmargin\labelsep
  \itemsep\z@\parsep\z@\topsep\z@\parskip\z@
  \usecounter{enumi}}
  \def\newblock{\hskip .11em plus .33em minus .07em}
  \sloppy\clubpenalty4000\widowpenalty4000
  \sfcode`\.=1000\relax}

\def\big#1{{\hbox{$\left#1\vcenter
  to1.428\ht\strutbox{}\right.\n@space$}}}
\def\Big#1{{\hbox{$\left#1\vcenter
  to2.142\ht\strutbox{}\right.\n@space$}}}
\def\bigg#1{{\hbox{$\left#1\vcenter
  to2.857\ht\strutbox{}\right.\n@space$}}}
\def\Bigg#1{{\hbox{$\left#1\vcenter
  to3.571\ht\strutbox{}\right.\n@space$}}}

\@addtoreset{equation}{section}
\def\theequation{\thesection.\arabic{equation}}

\usepackage{theorem}
\theorembodyfont{\itshape}
\theoremstyle{change}
\newtheorem{Definition}[equation]{Definition}
\newtheorem{Theorem}[equation]{Theorem}
\newtheorem{Proposition}[equation]{Proposition}
\newtheorem{Lemma}[equation]{Lemma}
\newtheorem{Corollary}[equation]{Corollary}

\theorembodyfont{\normalfont\rm}
\newtheorem{Remark}[equation]{Remark}

\newenvironment{Point}[1]%
  {\ifx*#1\let\pointlabel\relax\else\def\pointlabel{#1}\fi
   \refstepcounter{equation}\trivlist%
   \item[\hskip\labelsep\bf\theequation\space\pointlabel\space]%
   \ignorespaces\it
  }{\relax}

\let\DS\displaystyle

\let\SR\stackrel
\def\Prod{\displaystyle\prod}

\def\cf{{\it cf.\space}}

\def\){\big)}
\def\({\big(}
\let\>\rangle
\let\<\langle

\let\iso\cong

\let\0\varnothing              

\let\realb@r\bar
\let\bar\overline

\let\gedom\trianglerighteq
\let\gdom\vartriangleright

\def\Ann{\mathop{\rm Ann}\nolimits}
\def\End{\mathop{\rm End}\nolimits}
\def\Hom{\mathop{\rm Hom}\nolimits}

\def\Z{{\mathbb Z}}
\def\Q{{\mathbb Q}}

\let\To\longrightarrow
\def\map#1#2{\,{:}\,#1\!\longrightarrow\!#2}
\def\mapsto{\!\longmapsto\!}

\usepackage[vcentermath]{youngtab}
\def\tab(#1){\mbox{\small$\young(#1)$}\,}
\def\ydiag(#1){\mbox{\small$\yng(#1)$}\,}

\author{Richard Dipper and Andrew Mathas} 

\title{Morita equivalences of Ariki--Koike algebras} 
\address{\bigskip\sc 
  \begin{tabular}{ll} \makeatletter
    R.D.&Mathematisches Institut B, Universit\"at Stuttgart,\\ 
        &D-70550 Stuttgart. Deutschland.\\
        &{\tt rdipper@mathematik.uni-stuttgart.de}\\[2mm] 
    A.M.&School of Mathematics F07, University of Sydney,\\ 
        &Sydney N.S.W. 2006.  Australia.\\ 
        &{\tt mathas@maths.usyd.edu.au}\makeatother
\end{tabular} }
\AMSClassification{16G99, 20C20, 20G05} 
\Support{This paper is a contribution to the DFG projects ``Algorithmic
number theory and algebra'' (DI-531/3-2) and  
``Quantisierte Koordinatenringe und
verallgemeinerte $q$--Schur--Algebren'' (DI-531/5-1). 
Both authors acknowledge support
from DFG. The second author was supported in part by  the
Sonderforschungsbereich 343 in Bielefeld and the University of Sydney
research grant URG~U0628.}

\def\Q{{\bf Q}}
\def\H{\mathscr H}
\def\Sch{\mathscr S}
\def\Hb{\H_b\otimes\H_{n-b}}
\def\Sb{\Sch_b(\Gamma)\otimes\Sch_{n-b}(\Gamma)}
\def\HSn{\mathscr H(\Sym_n)} 
\let\len\ell
 
\let\phi\varphi
\def\phiST{\phi_{\S\T}}

\def\Nlam{N^\lambda}
\def\Nlambar{{\bar N}{}^\lambda}

\def\tlam{{\t^\lambda}}
\def\tmu{{\t^\mu}}

\def\phiST{\phi_{\S\T}}

\def\s{\mathfrak s}
\def\t{\mathfrak t}
\def\u{\mathfrak u}
\def\v{\mathfrak v}
\def\S{\mathtt S}
\def\T{\mathtt T}

\def\comp{\mathop{\rm comp}\nolimits}
\def\cont{\mathop{\rm cont}\nolimits}
\def\rad{\mathop{\rm rad}\nolimits}
\def\res{\mathop{\rm res}\nolimits}
\def\rest{{\downarrow}}
\def\Shape(#1){\mathop{\rm Shape}\nolimits(#1)}
\def\Sym{\mathfrak S}

\def\rtuple#1{({#1}^{(1)},\dots,{#1}^{(r)})}

\def\Std{\mathop{\rm Std}\nolimits}
\def\SStd(#1){{\mathcal T}_0(#1)}
\def\lSStd(#1){{\mathcal T}_l(#1)}

\def\D{\mathscr D}

\def\Mod{\operatorname{\bf Mod}}
\def\Modb{\Mod_{\H\!(b)}}

\def\Lam(#1:#2){\Lambda(#1\,{:}\,#2)}
\def\Lamp(#1:#2){\Lambda^+(#1\,{:}\,#2)}
\def\Gaml(#1:#2){\Gamma_l(#1\,{:}\,#2)}
\def\Gamr(#1:#2){\Gamma_r(#1\,{:}\,#2)}
\def\Gamlp(#1:#2){\Gamma_l^+(#1\,{:}\,#2)}
\def\Gamrp(#1:#2){\Gamma_r^+(#1\,{:}\,#2)}

\newcommand{\Lab}{\Lambda_b^+}
\newcommand{\barLab}{\bar{\Lambda}_b{\!}^+}

\newcommand{\Ho}{\rm{H}}
\newcommand{\hatHo}{\hat{\rm{H}}}

\newcommand{\Twbn}{T_{w_{n-b,b}}}

\def\Mb{M^{\omega_b}}
\def\Nb{{\bar N}{}^b}

\begin{document}
\def\@date{}\maketitle
\ifClassificationOrSupport\let\@makefnmark\relax
  \footnotetext{\kern-6mm{\it\@classification}\@support}
\fi

\centerline{\it Was lange w\"ahrt, wird endlich gut.}

\begin{abstract}
We prove that every Ariki--Koike algebra is Morita equivalent to a
direct sum of tensor products of smaller Ariki--Koike algebras which
have $q$--connected parameter sets. A similar result is proved for
the cyclotomic $q$--Schur algebras. Combining our results with work
of Ariki and Uglov, the decomposition numbers for the Ariki--Koike
algebras defined over fields of characteristic zero are now known in
principle.
\end{abstract}

\section{Introduction}

Let $R$ be a commutative ring with $1$ and let
$q,Q_1,\dots,Q_r$ be elements of $R$ with $q$ invertible.
Let $\Q=(Q_1,Q_2,\dots,Q_r)$. The {\sf Ariki--Koike
algebra} $\H_{q,\Q}(n)$ is the associative unital
$R$--algebra with generators $T_0,T_1,\dots,T_{n-1}$
subject to the following relations 
$$\begin{array}{rlll} 
  (T_0-Q_1)\dots(T_0-Q_r) &=&0 \\
  T_0T_1T_0T_1&=&T_1T_0T_1T_0\\
  (T_i+1)(T_i-q) &=&0&\text{for $1\le i\le n-1$}\\
  T_{i+1}T_iT_{i+1}&=&T_iT_{i+1}T_i&\text{for $1\le i\le n-2$}\\ 
  T_iT_j&=&T_jT_i&\text{for $0\le i<j-1\le n-2$.} 
\end{array}$$

The main result of this paper states that up to Morita
equivalence the Ariki--Koike algebras depend only on the
$q$--orbits of the parameters in $\Q$. More precisely we have 
the following.

\begin{Theorem} Suppose that
$\Q=\Q_1\coprod\Q_2\coprod\dots\coprod\Q_\kappa$ $($disjoint
union$)$ is a partitioning $\Pi$ of the parameter set $\Q$ such that
$$f_\Pi(q,\Q)=\prod_{1\le \alpha<\beta\le\kappa}
         \prod_{Q_i\in\Q_\alpha \atop Q_j\in\Q_\beta}\
         \prod_{-n<a<n} (q^aQ_i-Q_j)$$
is an invertible element of $R$. Then $\H_{q,\Q}(n)$ is Morita 
equivalent to the algebra
$$\H_{q,\Pi}(n)
   =\bigoplus_{n_1,\dots,n_\kappa\ge0 \atop n_1+\dots+n_\kappa=n}
     \H_{q,\Q_1}(n_1)\otimes\H_{q,\Q_2}(n_2)\otimes
        \dots\otimes \H_{q,Q_\kappa}(n_\kappa).$$ 
\label{Morita equivalence}
\end{Theorem}

Notice that the polynomial $f_\Pi(q,\Q)$ is invertible only if
whenever there exist~$i$ and~$j$ with $q^aQ_i=Q_j$ for some $a$ with
$-n<a<n$ then $Q_i, Q_j\in\Q_\alpha$ for some $\alpha$. The
appearance of these polynomials is not surprising because it was
shown by Ariki~\cite{Ariki:ss} that over a field $\H$ is semisimple
if and only if $$P_\H(q,\Q)=\prod_{1\le i<j\le r}\prod_{-n<a<n}
(q^aQ_i-Q_j)\cdot \prod_{k=1}^n (1+q+\dots+q^{k-1})\ne0.$$ The
polynomial $P_\H(q,\Q)$ is the analogue of the Poincar\'e polynomial
for $\H$.  The factors of $P_\H(q,\Q)$ in the right hand product
determine whether or not the subalgebra $\H_q(\Sym_n)$ is
semisimple.

Permuting the parameters $Q_1,\dots,Q_r$ does not affect $\H$ up to
isomorphism. In addition, if $\Q'=c\Q=(cQ_1,\dots,cQ_r)$, where $c$
is any invertible element of~$R$, then $\H_{q,\Q}(n)$ and
$\H_{q,\Q'}(n)$ are isomorphic (replace $T_0$ with $c^{-1}T_0$).  We
now rephrase Theorem~\ref{Morita equivalence} so that it says that up to
Morita equivalence $\H$ depends only on the $q$--orbits of
$Q_1,\dots,Q_r$.  To state this precisely, say that~$\Q$ is {\sf
$q$--connected} if $Q_i=q^{a_i}$ for some integer~$a_i$ for each $i$
(for the purposes of Corollary~\ref{main2} we could require that $|a_i|<n$,
for each $i$, but this is not so important). Normally,~$\Q$ will not
be $q$--connected; however, up to a permutation of the
$\Q_\alpha$'s, there is a unique partitioning
$\Q=\Q_1\coprod\dots\coprod\Q_\kappa$ such that for
$\alpha=1,\dots,\kappa$ there exists an element~$c_\alpha\in R$ such
that $Q_i\in\Q_\alpha$ if and only if $Q_i=c_\alpha q^{a_i}$ for
some integer $a_i$; so,  $\Q_\alpha=c_\alpha\Q_\alpha'$ for
each~$\alpha$.  By the above remarks, if $c_\alpha$ is invertible
then $\Q_\alpha'$ is $q$--connected and $\H_{q,\Q_\alpha}(\Sym_n)$
is isomorphic to $\H_{q,\Q_\alpha'}(\Sym_n)$; so 
Theorem~\ref{Morita equivalence} implies the following.

\begin{Corollary} Suppose that each $Q_i$ is invertible for $1\le i\le r$.
Then the Ariki--Koike algebra $\H_{q,\Q}(n)$ is Morita equivalent
to a direct sum of tensor products of Ariki--Koike algebras which
have $q$--connected parameter sets.
\label{main2} 
\end{Corollary}

The importance of this result stems from the work of Ariki~\cite{Ariki:can}
which showed that if $R$ is a field of characteristic zero, $q\ne1$ and the
parameter set $\Q$ is $q$--connected then the decomposition numbers of
Ariki--Koike algebras can be computed in terms of the canonical basis of an
associated integral highest weight module for an affine quantum group.  In
addition,  Ariki~\cite{A:class} and Ariki and the second
author~\cite{AM:simples} used the results of \cite{Ariki:can} to classify
the irreducible representations of the Ariki--Koike algebras without any
restrictions on the field,~$q$ or~$\Q$; this was done by first reducing to
the case of $q$--connected parameter sets.  Corollary~\ref{main2} explains the
reduction of~\cite{AM:simples} by showing that it comes from a Morita
equivalence. 

In fact, we can do more than this because Uglov~\cite{Uglov}
(extending the ideas of~\cite{LT:Schur2}), gave an algorithm for
computing the decomposition matrices of the Ariki--Koike algebras
which satisfy the restrictions of Ariki's paper~\cite{Ariki:can}. In
the spirit of the Kazhdan--Lusztig conjectures, Uglov's algorithm
involves computing certain affine parabolic Kazhdan--Lusztig
polynomials and evaluating them at~$1$.  Combining \cite{Uglov} with
Proposition~\ref{spechtmultplyup}(iii) below we obtain the following.

\begin{Corollary} Suppose that $R$ is a field of characteristic zero, $q\ne1$ and
$Q_i\ne0$ for $1\le i\le r$. Then the decomposition matrix of
$\H_{q,\Q}(n)$ is known.

\end{Corollary}

Next consider the case where there exists an integer $a$ with 
$Q_i\ne q^aQ_j$ if and only if $i=j$; then Theorem~\ref{Morita equivalence} says
that $\H_{q,\Q}(n)$ is Morita equivalent to a direct sum of tensor
products of Hecke algebras of type $A$. This special case is a result of
Du and Rui~\cite[Theorem 4.14]{DuRui:akmorita}. 

\begin{Corollary}[Du--Rui] \label{durui} Suppose that 
$\prod_{1\le i<j\le r}\prod_{-n<a<n} (q^aQ_i-Q_j)$
is an invertible element of $R$. Then $\H_{q,\Q}(n)$ is Morita
equivalent to
$$\bigoplus_{n_1,\dots,n_r\ge0 \atop n_1+\dots+n_r=n}
\H_q(\Sym_{n_1})\otimes\H_q(\Sym_{n_2})\otimes\dots
             \otimes\H_q(\Sym_{n_r}).$$ 

\end{Corollary}

In fact, when $r=2$ this is a theorem of James and the first named author
\cite[Theorem 4.17]{DJ:B}. Although this paper is largely motivated
by~\cite{DJ:B}, the techniques we use are very different.  As
in~\cite{DJ:B} we explicitly construct the projective generator of $\H$
which induces the Morita equivalence of the main theorem; however,
unlike~\cite{DJ:B,DuRui:akmorita} we do this by adapting the standard basis
of~$\H$ (from~\cite{DJM:cyc}) to give bases for the family of
projective modules which describe the projective generator. This yields
precise information, such as Specht series for these modules, which is new
even in the special cases considered previously~\cite{DJ:B,DuRui:akmorita}.
Another consequence is that we are able to extend our results to the
cyclotomic $q$--Schur algebras $\Sch_{q,\Q}(n)$ of \cite{DJM:cyc}. 

\begin{Theorem} Suppose that $\Q=\Q_1\coprod\Q_2\coprod\dots\coprod\Q_\kappa$ 
$($disjoint union$)$ is a partitioning $\Pi$ of the parameter set $\Q$ 
such that $f_\Pi(q,\Q)$ is an invertible element of $R$. 
Then $\Sch_{q,\Q}(n)$ is Morita equivalent to the algebra
$$\Sch_{q,\Pi}(n)
   =\bigoplus_{\SR{n_1,\dots,n_\kappa\ge0}{n_1+\dots+n_\kappa=n}}
      \Sch_{q,\Q_1}(n_1)\otimes\Sch_{q,\Q_2}(n_2)
         \otimes\dots\otimes \Sch_{q,\Q_\kappa}(n_\kappa).$$
\label{main3} 
\end{Theorem}

Actually, a slightly more general statement is possible; see
Theorem~\ref{genmain3}.

Again some special cases of Theorem~\ref{main3} were known previously.  First, if
\mbox{$r=2$} then $\H$ is an Iwahori--Hecke algebra of type $B$ and in this
case Theorem~\ref{main3} can be deduced for special parameter sets from the results
of \cite{DGr,GrHi}.  The case $r=2$ is important because it has
implications for representation theory of symplectic and unitary groups;
see~\cite{DF:I,DF:II}. Secondly, Ariki~\cite{Ariki:cycqschur} (see also Du
and Rui~\cite{DuRui:akmorita}), have proved Theorem~\ref{main3} under the
assumptions of Corollary~\ref{durui}; here the cyclotomic $q$--Schur algebra arises
as a quotient of a quantum group of type $A$ acting on ``$q$--tensor
space'' as the centralizing algebra of the Ariki--Koike algebra.  

In order to prove Theorem~\ref{Morita equivalence} we first observe that  it
is enough to prove the following much simpler, but equivalent,
result.

\begin{Theorem} Fix an integer $s$ with $1\le s\le r$ and suppose that 
$$f_s(q,\Q)=\prod_{1\le i\le s<j\le r}\ 
            \prod_{-n<a<n}(q^aQ_i-Q_j)$$
is an invertible element of $R$. Then $\H_{q,\Q}$ is Morita equivalent 
to 
$$\H_{q,s,\Q}(n)=\bigoplus_{b=0}^n\H_{q,(Q_{s+1},\dots,Q_r)}(\Sym_b) 
                 \otimes\H_{q,(Q_{1},\dots,Q_s)}(\Sym_{n-b}). $$
\label{d=2}
\end{Theorem}

The general case follows by iterating Theorem~\ref{d=2}, using the remarks before
Corollary~\ref{main2}.

The proof of Theorem~\ref{d=2} is based on an explicit decomposition of
$\H=\H_{q,\Q}(n)$ into a direct sum of projective right ideals $V^b$
using the standard basis of $\H$. The different projective modules
$V^b$ have no direct summand in common and $V=\bigoplus_{b=0}^n V^b$
is a progenerator for $\H$. In Theorem~\ref{projective decomposition} we
give an explicit formula for the multiplicity of each $V^b$ in the
regular representation of~$\H$  and in Theorem~\ref{endomorphism ring} we
show that the endomorphism ring of $V$ is $\H_{q,s,\Q}(n)$; this
proves Theorem~\ref{d=2} and that the Morita equivalences of Theorem~\ref{d=2} are
given by the functors $-\otimes_{\H_{q,s,\Q}}V$ and $\Hom_\H(V,-)$.
The explicit description of these functors enables us to extend our
results to the cyclotomic $q$-Schur algebras in section~5.  

\section{The standard basis theorem}

Let $\Sym_n$ be the symmetric group on $\{1,2,\dots,n\}$, acting
from the right, and let $s_1,\dots,s_{n-1}$ be the standard Coxeter
generators of $\Sym_n$; that is, $s_i=(i,i+1)$ for all~$i$. If
$w\in\Sym_n$ write $w=s_{i_1}s_{i_2}\dots s_{i_k}$ and say this
expression is {\sf reduced} if $k$ is minimal; in this case, $k$ is the
{\sf length} of $w$ and we write $\len(w)=k$ and define
$T_w=T_{i_1}T_{i_2}\dots T_{i_k}$.  Let $\HSn$ be the $R$--span of
$\set{T_w|w\in W}$; then $\HSn$ is a free subalgebra of $\H$ of rank
$n!$ which is isomorphic to the Iwahori--Hecke algebra of~$\Sym_n$.
The Iwahori--Hecke algebra~$\H(\Sym_n)$ is described in detail in
\cite{M:ULect}. 

Let $L_1=T_0$ and for $1\le i<n$ set $L_{i+1}=q^{-1}T_iL_iT_i$.  These
elements satisfy the following fundamental relations (see~\cite[(3.3)]{AK1}
and \cite[(2.1)]{DJM:cyc}).

\begin{Point}*
Suppose that $1\le i\le n-1$ and $1\le j\le n$. Then
\begin{enumerate}\item $L_i$ and $L_j$ commute.  
\item $T_i$ and $L_j$ commute
if $i\ne j-1, j$.  
\item $T_i$ commutes with $L_iL_{i+1}$ and
with $L_i+L_{i+1}$.  
\item If $a\in R$ and $i\ne j$ then $T_i$
commutes with $(L_1-a)(L_2-a)\dots(L_j-a)$.  
\end{enumerate}\label{L-comm}

\end{Point}

The importance of these elements derives from the following
result.

\begin{Point}{(Ariki--Koike~\cite[Theorem 3.10]{AK1})} The algebra $\H$ 
is free as an $R$--module with basis 
$$
\set{L_1^{d_1}L_2^{d_2}\dots L_n^{d_n}T_w|
w\in\Sym_n \text{and} 0\le d_m\le r-1\text{for} m=1,2,\dots,n}.
$$ 
In particular, $\H$ is free of rank $r^nn!$
\label{AK basis} 
\end{Point}

Note that because $q$ is invertible so are the elements $T_i$, 
for $i=1,\dots,n-1$; explicitly,
$T_i^{-1}=q^{-1}(T_i-q+1)$. Consequently, $T_w$ is invertible for
all $w\in\Sym_n$.

Let $*\map\H\H$ be the anti--automorphism of $\H$
determined by $T_i^*=T_i$ for $i=0,1,\dots,n-1$. Then
$T_w^*=T_{w^{-1}}$ for all $w\in\Sym_n$ and $L_i^*=L_i$
for $i=1,2,\dots,n$.

As we next recall, the Ariki--Koike algebra has another basis which is
better adapted to the study of its representation theory; this basis
is {\sf cellular}, in the sense of Graham and Lehrer~\cite{GL}.
(Graham and Lehrer were the first to construct a cellular basis of
$\H$; the basis we use is due to Gordon James and the
authors~\cite{DJM:cyc}.)

A {\sf composition} of an integer $m\ge0$ is an ordered sequence of
non--negative integers $\tau=(\tau_1,\tau_2\dots)$ such that
$|\tau|\:=\sum_{i\ge1} \tau_i=m$; if the sequence is non--increasing then
$\tau$ is a {\sf partition} of $m$. A {\sf multicomposition} of $n$ (with
$r$--components) is an ordered $r$--tuple $\mu=(\mu^{(1)},\dots,\mu^{(r)})$
of compositions such that $|\mu^{(1)}|+\dots+|\mu^{(r)}|=n$. If each
$\mu^{(i)}$ is a partition then $\mu$ a {\sf multipartition}.  Let
$\Lambda=\Lam(n:r)$ be the set of multicompositions of $n$ with
$r$--components and let~$\Lambda^+=\Lamp(n:r)\subset\Lambda$ be the set of
multipartitions of $n$ with $r$--components.

The {\sf diagram} of a multicomposition $\mu$ is the set
$$[\mu]=\set{(i,j,k)|1\le k\le r\And i\ge1\And 1\le
                        j\le\mu^{(k)}_i},$$ 
which we think of as an ordered $r$--tuple of boxes in the
plane. For example, if $\mu=((3,1),(1^2),(2,1))$ then
$$[\mu]=\Bigg(\,\ydiag(3,1),\ \ydiag(1,1),\ \ydiag(2,1)\,\Bigg).$$ 
In this way, we talk of the rows and columns (of the components) of
$\mu$.

Given two multicompositions $\lambda$ and $\mu$ say that
$\lambda$ {\sf dominates} $\mu$, and write
$\lambda\gedom\mu$, if for $1\le c\le r$ and for all~$i\ge 1$
$$\sum_{b=1}^{c-1}|\lambda^{(b)}|
      +\sum_{j=1}^i\lambda^{(c)}_j\ge
  \sum_{b=1}^{c-1}|\mu^{(b)}|+\sum_{j=1}^i\mu^{(c)}_j.$$ 
If $\lambda\gedom\mu$ and $\lambda\ne\mu$ we write
$\lambda\gdom\mu$. This defines a partial order on the sets of
multicompositions and multipartitions of $n$.

If $\mu$ is a multicomposition of $n$ then a {\sf $\mu$--tableau} is a
map $\t\map{[\mu]}\{1,2\dots,n\}$; we write $\Shape(\t)=\mu$.
Generally, we shall think of tableaux as labelled diagrams; for
example, when $\mu=((3,1),(1^2),(2,1))$ three $\mu$--tableaux are
$$\Big(\ \tab(123,4),\,\tab(5,6),\,\tab(78,9)\ \Big),\quad 
  \Big(\ \tab(358,4),\,\tab(2,9),\,\tab(16,7)\ \Big)\quad\And\quad 
  \Big(\ \tab(358,4),\,\tab(7,6),\,\tab(12,9)\ \Big).$$ 
Given a tableau $\t$ we will also write $\t=\rtuple\t$ and call
$\t^{(i)}$ the $i$th {\sf component} of $\t$. Similarly, given an
integer $m$ with $1\le m\le n$ we write $\comp_\t(m)=i$ if $m$ appears
in the $i$th component $\t^{(i)}$ of $\t$. If $\s$ is another tableau
we write $\comp_\t=\comp_\s$ if $\comp_\t(m)=\comp_\s(m)$ for all $m$
with $1\le m\le n$.

A $\mu$--tableau $\t$ is {\sf standard} if the entries in each row
of each component of~$\t$ increase from left to right and the entries
in each column of each component of $\t$ increase from top to bottom
(of the tableaux above only the first two are standard). For each
multipartition $\lambda$ let $\Std(\lambda)$ be the set of standard
$\lambda$--tableaux.

Suppose $\t$ is a standard $\lambda$--tableaux and $\s$ a standard
$\mu$--tableau for $\lambda, \mu\in\Lambda$. Given an integer $m$
with $1\le m\le n$ let $\t\rest m$ be the subtableau of $\t$ which
contains the entries $1,2,\dots,m$ and similarly for $\s\rest m$.
Then $\s\gedom\t$ if $\Shape(\s\rest m)\gedom\Shape(\t\rest m)$, for
$m=1,2,\dots,n$, and we say that $\s$ {\sf dominates} $\t$. Again we
write $\s\gdom\t$ if $\s\gedom\t$ and $\s\ne\t$.

Let $\tmu$ be the unique $\mu$--tableau such that $\tmu\gedom\t$ for
all $\mu$--tableau $\t$. Then~$\tmu$ is the tableau which has the
numbers $1,2,\dots,n$ entered in order along the rows of $[\mu]$. Note
that the symmetric group $\Sym_n$ acts from the right on the set of
$\mu$--tableaux. Let $\Sym_\mu$ be the row stabilizer of the tableau
$\t^\mu$; then $\Sym_\mu$ is a parabolic subgroup of $\Sym_n$. In the
example above, where $\mu=((3,1),(1^2),(2,1))$, the
first of the tableaux listed is~$\tmu$ and  
$\Sym_\mu=\Sym_3\times\Sym_1\times\Sym_1\times\Sym_1
             \times\Sym_2\times\Sym_1\hookrightarrow\Sym_9$ 
(obvious embedding); we will always identify $\Sym_\mu$ with a
subgroup of $\Sym_n$ in this way.

For each $\mu$--tableau $\t$ let $d(\t)$ be the unique element of
$\Sym_n$ such that $\t=\tmu d(\t)$. Then, by \cite[Lemma 1.4]{DJ:reps},
$d(\t)$ is a distinguished right coset representative of~$\Sym_\mu$
in $\Sym_n$; that is, $\len(wd(\t))=\len(w)+\len(d(\t))$ for all
$w\in\Sym_\mu$.
 
Given an $r$--tuple $\mathbf a=(a_1,a_2,\dots, a_r)$ of
integers, with $0\leq a_i\leq n$ for all $i$, let
$u_{\mathbf a}=u_{a_1,1}u_{a_2,2}\dots u_{a_r,r}$ 
where $u_{a,t}=\prod_{k=1}^{a}(L_k-Q_t)$ for any $a$ and $t$.
If $\mu$ is a multicomposition of $n$ let $u^+_\mu=u_{\mathbf a}$
where $\mathbf a=(a_1,a_2,\dots, a_r)$ is the sequence with
$a_t=|\mu^{(1)}|+\dots+|\mu^{(t-1)}|$ for $t=1,2,\dots,r$.

Suppose that $\lambda$ is a multipartition of $n$ and let
$x_\lambda=\sum_{w\in\Sym_\lambda}T_w$ and set
$m_\lambda=u^+_\lambda x_\lambda$. If $\s$ and $\t$ are standard
$\lambda$--tableaux define $m_{\s\t}=T_{d(\s)}^*m_\lambda
T_{d(\t)}^{\phantom*}$. Then we have the following.

\begin{Point}{(Dipper--James--Mathas \cite[Theorem 3.26]{DJM:cyc})}  
Let 
$$ \mathcal M=\set{m_{\s\t}|\s,\t\in\Std(\lambda)\ForSome
                            \lambda\in\Lambda^+}.
$$
Then $\mathcal M$ is a cellular basis of $\H$.
\label{standard basis}
\end{Point}

The basis $\mathcal M$ is called the {\sf standard basis} of $\H$.

As in \cite{DJM:cyc}, let $\Nlam$ be the $R$--module with basis the
set of all $m_{\s\t}\in\mathcal M$ where $(\s,\t)$ runs over all
pairs of standard $\mu$--tableaux with $\mu\gedom\lambda$;
similarly, let $\Nlambar$ be the $R$--module with basis the set of
$m_{\u\v}$ where $\u$ and $\v$ are standard $\mu$--tableau and 
$\mu\gdom\lambda$.  From the theory of cellular algebras
we obtain the following corollary of (\ref{standard basis}).

\begin{Point}{(\cite[Corollary 3.22]{DJM:cyc})} The $R$--modules $\Nlam$ 
and $\Nlambar$ are two--sided ideals of~$\H$.  
\label{idealsequence}
\end{Point}

Given a multipartition $\lambda$ the {\sf Specht module} $S^\lambda$
is the submodule of $\H/\Nlambar$ defined by $S^\lambda=z_\lambda\H$
where $z_\lambda=m_\lambda+\Nlambar$. The theory of cellular
algebras~\cite{GL,M:ULect} shows that $S^\lambda$ is free of rank
$|\Std(\lambda)|$ (with basis
$\set{m_{\tlam\t}+\Nlambar|\t\in\Std(\lambda)}$), and that there is
an intrinsically defined symmetric $\H$--invariant bilinear form 
$(\ ,\ )$ on~$S^\lambda$; that is $(uh,v)=(u,vh^*)$ for $u,v\in
S^\lambda$ and $h\in\H$. Let $\rad S^\lambda$ be the radical of this
form and set
$D^\lambda=S^\lambda/\rad S^\lambda$.  Then~$D^\lambda$ is an
$\H$--module; moreover, the following is true.

\begin{Point}{\cite{GL,DJM:cyc}} Suppose that $R$ is a field. 
\begin{enumerate}\item If $\lambda$ is a multipartition of $n$ then $D^\lambda$ is
either $(0)$ or absolutely irreducible.
\item $\set{D^\lambda|\lambda\in\Lambda^+\And D^\lambda\ne (0)}$
is a complete set of pairwise non--isomorphic irreducible $\H$--modules.
\end{enumerate}

\end{Point}

Given multipartitions $\lambda$ and $\mu$ with $D^\mu\ne(0)$ let
$d_{\lambda\mu}=[S^\lambda:D^\mu]$ be the decomposition multiplicity
of the simple module $D^\mu$ in the Specht module $S^\lambda$. The
matrix $(d_{\lambda\mu})$ is the {\sf decomposition matrix} of~$\H$.
Importantly, the decomposition matrix of $\H$ is unitriangular; more
precisely, we have the following.

\begin{Point}{\cite{GL,DJM:cyc}}
Suppose that $R$ is a field and that $\lambda$ and $\mu$ are
multipartitions such that $D^\mu\ne(0)$. Then $d_{\mu\mu}=1$ and
$d_{\lambda\mu}\ne0$ only if $\lambda\gedom\mu$.
\label{unitriangular}
\end{Point}

Let $\s$ be a tableau and suppose that $1\le k\le n$ appears in row
$i$ and column~$j$ of the $c$th component $\s^{(c)}$ of $\s$. Then
the {\sf residue} of $k$ in $\s$ is $\res_\s(k)=q^{j-i}Q_c$. The
following result underpins much of what follows.

\begin{Point}{\cite[Prop.~3.7]{JM:cyc-Schaper}} Let $\s$ and $\t$ be
standard $\lambda$--tableaux, where $\lambda\in\Lambda^+$, and
suppose that $1\leq k\leq n$. Then there exist $a_\u\in R$ such that
$$
L_km_{\s\t}\equiv \res_\s(k)m_{\s\t}
               +\sum_{\u\gdom\s}a_\u m_{\u\t} \mod\Nlambar.
$$ 
\label{L_k action}
\end{Point}

For each multipartition $\lambda$ define its {\sf content} to be
the sequence $\cont(\lambda)=(c_r)_{r\in R}$ where $c_r$ is the
number of nodes $x$ of the diagram $[\lambda]$ with $\res(x)=r$.
By \cite[Theorem~3.7(ii)]{GL} all of the irreducible constituents of
$S^\lambda$ belong to the same block. Furthermore, by
(\ref{L-comm}) every symmetric polynomial  $f(L)$, in
$L_1,L_2,\dots,L_n$ belongs to the centre of $\H$; therefore, by
Schur's Lemma, $f(L)$ acts on $S^\lambda$ as multiplication by a
scalar, say $\alpha_f\in R$. By (\ref{L_k action}),
$f(L)m_\lambda\equiv\alpha_fm_\lambda\mod\Nlambar$; so it
follows that $\alpha_f$ depends only upon the content of
$\lambda$. Hence, we have the following result.

\begin{Corollary}[\protect{Graham--Lehrer~\cite[Prop.~5.9(ii)]{GL}}] Suppose that
$\lambda$ and $\mu$ are multipartitions of $n$. Then $S^\lambda$
and $S^\mu$ belong to the same block only if
$\cont(\lambda)=\cont(\mu)$.
\label{content}
\end{Corollary}

Grojnowski~\cite{Groj:AKblocks} has recently shown that $S^\lambda$
and $S^\mu$ are in the same block if and only if
$\cont(\lambda)=\cont(\mu)$. (The definition of residue must be
modified slightly in the case $q=1$.) When $r=1$ this result was
already known by \cite{DJ:blocks,JM:Schaper}.

For each multicomposition $\mu$ of $n$ let
$M^\mu=m_\mu\H$.  To describe how the standard basis of
(\ref{standard basis}) can be modified to give a basis of $M^\mu$ we
need to generalize the notion of tableau. Suppose that $\lambda$
is a multipartition and $\mu$ is a multicomposition of $n$. A 
{\sf $\lambda$--tableau of type $\mu$} is a map
$\S\map{[\lambda]}\{1,\dots,n\}\times\{1,\dots,r\}$ such that, for
all~$(i,k)$, $\mu_i^{(k)}=\#\set{x\in[\lambda]|\S(x)=(i,k)}$; as
before, we will think of $\S=\rtuple\S$ as a labelling of
$[\lambda]$ with ordered pairs of integers $(i,k)$. A
$\lambda$--tableau $\S$ of type~$\mu$ is {\sf semistandard} if for
$c=1\dots,r$ the entries in the $c$th component $\S^{(c)}$ of $\S$
are (i)~non--decreasing along the rows; (ii)~strictly increasing
down the columns and, (iii)~no entry in $\S^{(c)}$ has the form
$(i,k)$ with $k<c$. Let $\SStd(\lambda,\mu)$ be the set of
semistandard $\lambda$--tableaux of type $\mu$.

Given a standard $\lambda$--tableau $\s\map{[\lambda]}\{1,\dots,n\}$ 
define $\mu(\s)$ to be the  $\lambda$--tableau of type $\mu$
obtained from $\s$ by replacing each entry $m$ in $\s$
by $(i,k)$, if~$m$ appears in row $i$ of the $k$th component of
$\t^\mu$. Observe that in general the entries in each column of 
$\mu(\s)$ will only be non--decreasing so that $\mu(\s)$ need not be
semistandard.

\begin{Point}rm For example, if we let $\omega=\((0),\dots,(0),(1^n)\)$ then the
tableau $\omega(\s)$ has entries of the form $(i,r)$ where 
$1\le i\le n$. Consequently, $\s\mapsto\omega(\s)$ gives a bijection
from the set of standard $\lambda$--tableaux to the set of
semistandard $\lambda$--tableaux of type $\omega$. Henceforth, we
identify $\lambda$--tableau of type $\omega$ with the
$\lambda$--tableaux that are maps from $[\lambda]$ to
$\{1,\dots,n\}$; we also use lower case letters $\s,\t,\dots$ to
denote $\lambda$--tableaux (of type $\omega$) and upper case letters
$\S,\T,\dots$ for tableaux of arbitrary type.
\label{type omega}

Before we can  state the basis theorem for $M^\mu$ we need one more
definition. Let $\S$ be a semistandard tableau of type $\mu$ and a
$\lambda$--tableau $\t$; set
$$m_{\S\t}=\sum_{\SR{\s\in\Std(\lambda)}{\mu(\s)=\S}} m_{\s\t}.$$

\end{Point}

\begin{Point}{\cite[Theorem 4.14]{DJM:cyc}}
Suppose that $\mu$ is a multicomposition of $n$. Then $M^\mu$ is
free as an $R$--module with basis 
$$\set{m_{\S\t}|\S\in\SStd(\lambda,\mu),
               \t\in\Std(\lambda)\ForSome\lambda\in\Lambda^+}.$$
\label{semistandard basis of Mmu}
\end{Point}

The last result that we shall need gives a basis for
$\Hom_\H(M^\nu,M^\mu)$. Given a semistandard $\lambda$--tableau $\S$
of type $\mu$ and a semistandard tableau $\T$ of type $\nu$ let
$$m_{\S\T}
  =\sum_{\SR{\s,\t\in\Std(\lambda)}{\mu(\s)=\S,
  \nu(\t)=\T}}m_{\s\t}$$
and define $\phiST\map{M^\nu}M^\mu$ to be the $\H$--module
homomorphism given by
$\phiST(m_\nu h)=m_{\S\T}h$ for all $h\in\H$.
It is not completely obvious that $\phiST$ even belongs to
$\Hom_\H(M^\nu,M^\nu)$; nonetheless, the following is true.
 
\begin{Point}{\cite[Theorem 6.6(i)]{DJM:cyc}}
Suppose that $\mu$ and $\nu$ are multicompositions of $n$. Then
$\Hom_\H(M^\nu,M^\nu)$ is free as an $R$--module with basis
$$\set{\phiST|\S\in\SStd(\lambda,\mu)\And\T\in\SStd(\lambda,\nu)
                  \ForSome \lambda\in\Lambda^+}.$$
\label{SStd basis}
\end{Point}

\section{A projective generator for $\H$}

For the remainder of the paper we fix an integer $s$ with $1\le s\le r$.

Following \cite{DJ:B}, given integers $i$ and $j$ with 
$1\le i<j<n$ define 
$$ s_{i,j}=s_is_{i+1}\dots s_{j-1} \quad\And\quad
s_{j,i}=s_{i,j}^{-1}=s_{j-1}s_{j-2}\dots s_i.
$$
Note that $s_{j,i}$ is the cycle $(i,i+1,\dots,j)$.
We abbreviate the corresponding elements of $\H$ by 
$T_{i,j}=T_{s_{i,j}}$ and $T_{j,i}=T_{s_{j,i}}$; so,
$T_{j,i}=T_{i,j}^*$. In passing, we remark that
$L_i=q^{1-i}T_{i,1}T_0T_{1,i}$ for $i=1,\dots,n$.

If $a$ and $b$ are non--negative integers we define
$w_{a,b}=(s_{a+b,1})^b$; in particular, $w_{a,0}=w_{0,b}=1$.
Written as a permutation,
$$w_{a,b}=\left(\begin{array}{cccccccc} 
               1 &2 &\dots&a&a+1&a+2&\dots&a+b\\ 
              b+1&b+2 &\dots&a+b&1&2&\dots&b
        \end{array}\right).$$
For later use we note that
$w_{a,b}^{-1}=w_{b,a}$; consequently, $T_{w_{a,b}}^*=T_{w_{b,a}}$.

\begin{Lemma} Let $a$ and $b$ be non--negative integers
with $0\le a+b\le n$ and suppose that $i$ is an integer such
that $i\ne a$ and $1\le i<a+b$. Then
$$ T_iT_{w_{a,b}} =T_{w_{a,b}}T_{(i)w_{a,b}}
            = \begin{cases}
                T_{w_{a,b}}T_{i+b},&\If 1\le i<a,\\
                T_{w_{a,b}}T_{i-a},&\If a<i<a+b.
              \end{cases}
$$
\label{w_b properties}

\end{Lemma}

\proof For any $w\in\Sym_n$ we have 
$s_iw=ww^{-1}s_iw=w(iw,(i+1)w)$; setting $w=w_{a,b}$ we find
$$
s_iw_{a,b}=\begin{cases} w_{a,b}s_{b+i},&\If 1\leq i< a,\\  
                         w_{a,b}s_{i-a},&\If a<i<a+b.
\end{cases}
$$
The result follows by observing that
$\len(s_iw_{a,b})=\len(w_{a,b})+1$ since $w_{a,b}$ is a
distinguished right coset representative of $\Sym_{(a,b)}$ in
$\Sym_n$; see~\cite[2.7]{DJ:B}.
\endproof

Equation \cite[2.9]{DJ:B} provides a recursive way to find a 
reduced expression for the permutation $w_{a,b}$. As a consequence
we obtain the following Lemma.

\begin{Lemma} Suppose that $1\leq b\leq n-1$. Then there exists
$\tilde{w}\in\Sym_{(1,n-1)}$ such that 
$w_{n-b,b}=\tilde{w}s_{1,b+1}$
and $\len(w_{n-b,b})=\len(\tilde{w})+\len(s_{1,b+1})$. In
particular, $T_{w_{n-b,b}}= T_{\tilde{w}}T_{1,b+1}$.
\label{w_b properties 2}

\end{Lemma}

\proof By \cite[2.9]{DJ:B} if $a\ge0$ and $0\le a+b<n$
then $w_{a,b}=s_{a,a+b}w_{a-1,b}$ and 
$\len(w_{a,b})=\len(s_{a,a+b})+\len(w_{a-1,b})$.
Therefore, 
\begin{align*}
w_{n-b,b}&=s_{n-b,n}w_{n-b-1,b}=s_{n-b,n}s_{n-b-1,n-1}w_{n-b-2,b}
          =\dots\\
         &=s_{n-b,n}s_{n-b-1,n-1}\dots s_{2,b+2}w_{1,b}
          =(s_{n-b,n}\dots s_{2,b+2})s_{1,b+1},
\end{align*}
with the lengths adding throughout. The Lemma follows.
\endproof

For the next definition recall that we have fixed an integer $s$
with $1\le s\le r$. For each integer $b$ with $0\le b\le n$ define
\begin{align*}
u_{n-b}^-&=\prod_{t=1}^{s}(L_1-Q_t)(L_2-Q_t)\dots(L_{n-b}-Q_t)\\ 
\intertext{and}
u_{b}^+&=\prod_{t=s+1}^{r}(L_1-Q_t)(L_2-Q_t)
              \dots(L_{b}-Q_t).
\end{align*}
These elements are special instances of the elements $u_{\mathbf a}$
introduced in the previous section.

We now define the modules which are the cornerstone upon which
this paper is built.

\begin{Definition} Suppose that $0\le b\le n$ and define
$v_b=u_{n-b}^- \Twbn  u_b^+$ and let $V^b=v_b\H$.

\end{Definition}

Ultimately we shall show that $V^b$ is a projective $\H$--module and
that its endomorphism ring is isomorphic to a tensor product of
smaller Ariki--Koike algebras; this will imply our main result.  At
this point it is not even clear that~$v_b$ is non--zero; we will
deduce this important fact latter. We begin by establishing some key
properties of $v_b$.

\begin{Proposition} Suppose that $0\le b\le n$. 
\begin{enumerate}\item If $1\le i<n-b$ then $T_iv_b=v_b T_{i+b}$.
\item If $n-b<i\le n$ then $T_iv_b=v_b T_{i-n+b}$.
\item If $1\le k\le n-b$ then $L_k v_b=v_bL_{k+b}$.
\item If $n-b+1\le k\le n$ then $L_k v_b=v_bL_{k-n+b}$.
\end{enumerate}\label{v_b properties}

\end{Proposition}

\proof First, observe that parts (i) and (ii) follow from parts 
(i) and (ii) of Lemma~\ref{w_b properties}, respectively, together
with~(\ref{L-comm})(iv). 

Next, consider part (iii). If $b=0$ or $b=n$ then $v_b$ is central in $\H$
by parts (ii) and (iii) of (\ref{L-comm}); so our claims follow in these two
cases and we may assume that $1\leq b\leq n-1$. 

We first consider the case $k=1$; that is, $L_1v_b$. We write
$w_{n-b,b}=\tilde{w}s_{1,b+1}$, as in Lemma~\ref{w_b properties 2}; then
$\tilde{w}\in S_{(1,n-1)}$ so that $L_1$ and $T_{\tilde{w}}$ commute
by (\ref{L-comm}). Now each $L_i$ commutes with $u_{n-b}^-$ since the
$L_i$ generate an abelian subalgebra of $\H$; therefore,
\begin{align*}
L_1v_b &= L_1u_{n-b}^-\Twbn u_b^+\\
&= u_{n-b}^-L_1T_{\tilde{w}}T_{1,b+1}u_b^+\\
&= u_{n-b}^-T_{\tilde{w}}L_1T_{1,b+1}u_b^+\\
&=u_{n-b}^-T_{\tilde{w}}T_{b+1,1}^{-1}T_{b+1,1}L_1T_{1,b+1}u_b^+\\ 
&= q^{b}u_{n-b}^-T_{\tilde{w}}T_{b+1,1}^{-1}L_{b+1}u_b^+\\
&= q^{b}u_{n-b}^-T_{\tilde{w}}(T_1^{-1}
         \dots T_{b}^{-1})u_b^+L_{b+1}.
\end{align*}
Now $q^b(T_1^{-1}\dots T_{b}^{-1})=
(T_1-q+1)\dots(T_{b}-q+1)=T_1\dots T_{b}+h$, where $h\in\H$ is
an $R$--linear combination of terms of the form $T_xT_y$ such that
$(x,y)$ is an element of $\Sym_i\times\Sym_{b-i}=\Sym_{(i,b-i)}$
for some $i$ with $0<i<b$. Write $u_{n-b}^-=u_1\tilde{u}$, where
$u_1=\prod_{t=1}^s(L_1-Q_t)$, and suppose $(x,y)\in\Sym_{(i,b-i)}$
for some~$i>0$. Then $T_xT_y=T_{xy}=T_{yx}=T_yT_x$ and, by
(\ref{L-comm}), $T_x$ commutes with $u_b^+$, $T_y$ commutes with
$\tilde u$ and $T_{\tilde w}$ commutes with $u_1$. Therefore,
$$u_{n-b}^-T_{\tilde w}T_{xy}u_b^+
    =u_1\tilde{u}T_{\tilde w}T_yT_xu_b^+
    =\tilde{u}T_{\tilde w}u_1T_yT_xu_b^+
    =\tilde{u}T_{\tilde w}T_yu_1u_b^+T_x
    =0;$$
the last equality following because $\prod_{t=1}^r(L_1-Q_t)=0$ is a
factor of $u_1u_b^+$.  Consequently,
$q^{b}u_{n-b}^-T_{\tilde{w}}hu_b^+L_{b+1}=0$.  Hence,
$$L_1v_b=u_{n-b}^-T_{\tilde{w}}(T_1\dots T_{b})u_b^+L_{b+1}
        =u_{n-b}^-\Twbn u_b^+L_{b+1}
        =v_bL_{b+1} 
$$
as claimed.

Next consider $L_kv_b$ for some $k$ with $1<k\le n-b$. 
By induction and part~(i),
$$L_kv_b =q^{-1}T_{k-1}L_{k-1}T_{k-1}v_b
         =q^{-1}T_{k-1}v_bL_{b+k-1}T_{b+k-1}
         =v_bL_{b+k},$$
proving (iii).

Part (iv) can be proved similarly; however, here is a better
argument. As there is no essential difference between $u_b^+$ and
$u_{n-b}^-$ --- and hence between $v_b$ and $v_b^*$ --- it follows
from part~(iii) that if $n-b<k\le n$ then 
$L_kv_b=(v_b^*L_k)^*=(L_{k-n+b}v_b^*)^*=v_bL_{k-n+b}$, giving (iv). 
\endproof

Define $\mathsf L_A(m)=\set{x\in A|xm=0}$ to be the {\sf left annihilator}
of $m\in M$ in~$A$; $\mathsf L_A(h)$ is a left ideal of $A$. We also let
$\Ann_A(M)=\bigcap_{m\in M}\mathsf L_A(m)$ be the {\sf annihilator} of $M$;
this is an ideal of $A$. We will apply these definitions in the case where
$M=V^b$ and $A=H_b$ is the subalgebra of $\H$ generated by
$\set{T_i,T_j,L_k| 1\le i<n-b< j<n\And1\le k\le n}$.  Observe that by the
Proposition~$V^b$ is invariant under left multiplication by $H_b$ and
hence a left $H_b$-module.

\begin{Corollary} Suppose that $0\leq b\leq n$ and let the subalgebra $H_b$ of $\H$ 
be defined as above. Then $H_b$ acts on $V^b$ by left multiplication and
$\mathsf L_\H(v_b)\cap H_b=\Ann_{H_b}(V^b)$; consequently, the algebra 
$\hat H_b=H_b/\Ann_{H_b}(V^b)$ is a subalgebra of $\End_\H(V^b)$.
\label{endhom}
\end{Corollary}

In fact, in Theorem~\ref{endomorphism ring} below we will prove that 
$$\hat H_b\cong\End_\H(V^b)\cong
  \H_{q,(Q_1,\dots,Q_s)}(b)\otimes\H_{q,(Q_{s+1},\dots,Q_{r})}(n-b).$$
The next result shows that when $L_{n-b+1}$ and $L_1$act on $V^b$ by
left multiplication they each satisfy one of the relations of the
generators $1\otimes T_0$ and $T_0\otimes 1$, respectively, in the
tensor product above.  

\begin{Corollary} Suppose that $0\le b\le n$. Then
\begin{enumerate}\item $(L_1-Q_{s+1})\dots(L_1-Q_r)v_b=0$; 
\item $(L_{n-b+1}-Q_1)\dots(L_{n-b+1}-Q_s)v_b=0$;
\item $v_b(L_1-Q_1)\dots(L_1-Q_s)=0$; and,
\item $v_b(L_{b+1}-Q_{s+1})\dots(L_{b+1}-Q_r)=0$.
\end{enumerate}\label{v_b annihilators}

\end{Corollary}

\proof Parts (i) and (iii) follow from the relation
$\prod_{t=1}^r(L_1-Q_t)=0$ and the definition of $v_b$; for parts
(ii) and (iv) apply the Proposition~\ref{v_b properties} to parts (iii) and~(i)
respectively.
\endproof          

\begin{Corollary} Let  $0\le b < c\le n$. Then $u_{n-b}^-\Twbn u_c^+ = 0$.
\label{ubTucnull}

\end{Corollary}

\proof Let $h=\prod_{t=s+1}^r(L_{b+2}-Q_t)\cdots (L_c-Q_t)$. Then
\begin{align*}
u_{n-b}^-\Twbn u_c^+
   &= u_{n-b}^-\Twbn \prod_{t=s+1}^r(L_1-Q_t)(L_2-Q_t)\cdots (L_c-Q_t)\\
      &= u_{n-b}^-\Twbn u_b^+(L_{b+1}-Q_{s+1})\cdots (L_{b+1}-Q_r)h\\
      &= v_b(L_{b+1}-Q_{s+1})\cdots (L_{b+1}-Q_r)h=0,
\end{align*}
where the last equality comes from Corollary~\ref{v_b annihilators}(iv).
\endproof

We will study the ideals $V^b$ by thinking of them as quotients of
one of the modules $M^\mu$. Let
$\omega_b=(\omega_b^{(1)},\dots,\omega_b^{(r)})$ be the
multipartition of $n$ with 
$$\omega_b^{(t)}=\cases (1^b),     &\If t=s,\\
                        (1^{n-b}), &\If t=r,\\
                        (0),       &\Otherwise.
                     \endcases$$
Then $u_b^+=u^+_{\omega_b}=m_{\omega_b}$; consequently,
$v_b=u_{n-b}^-\Twbn m_{\omega_b}$ and $V^b$ is a quotient of
$\Mb$. This motivates the following definition.

\begin{Definition} Suppose that $0\le b\le n$. Let 
$\theta_b\map{\Mb}V^b$ be the map given by
$\theta_b(h)=u_{n-b}^-\Twbn h$ for all $h\in\Mb$.

\end{Definition}

Thus, $\theta_b$ is a surjective $\H$--module homomorphism from
$\Mb$ onto $V^b$. The map $\theta_b$ is the main tool we need to
understand the modules $V^b$; first we set up some notation.

Let $\lambda$ be a multipartition of $n$ and set
$$\Std_b(\lambda)
  =\set{\t\in\Std(\lambda)|\comp_\t(k)\le s\Whenever 1\le k\le b};$$ 
that is, $\t\in\Std_b(\lambda)$ if and only if the numbers
$1,2,\dots,b$ all appear in one of the first~$s$ components of $\t$.
Similarly, let
$$\Std_{b,n-b}(\lambda)=
   \set{\t\in\Std_b(\lambda)|\comp_\t(k)>s\Whenever b<k\le n}.$$
Let $\Lab=\set{\lambda|\lambda\in\Lambda^+\And
        |\lambda^{(1)}|+\dots+|\lambda^{(s)}|=b};$
then $\Std_{b,n-b}(\lambda)$ is non--empty if and only if
$\lambda\in\Lab$. On the other hand, $\Std_b(\lambda)$ is
non--empty if and only if 
$|\lambda^{(1)}|+\dots+|\lambda^{(s)}|\ge b$; we set
$$\barLab=\set{\lambda|\lambda\in\Lambda^+\And
        |\lambda^{(1)}|+\dots+|\lambda^{(s)}|>b}.$$
Then $\barLab$ is a coideal in $\Lambda^+$; that is, if
$\mu\in\Lambda^+$ and $\mu\gedom\lambda$ for some $\lambda\in\Lab$
then $\mu\in\Lab$. In contrast, $\Lab$ is not a coideal; however,
$\Lab\cup\barLab$ is a coideal and $\Lab$ and
$(\Lab\cup\barLab)/\barLab$ are isomorphic posets.

Observe that if $\lambda\in\Lab$ then
$\Std_{b,n-b}(\lambda)=\Std_b(\lambda)$; however, we will
continue to write $\Std_{b,n-b}(\lambda)$ in order to emphasize
the restrictions on the components of these tableaux. 

\begin{Lemma} Suppose that $0\le b\le n$. Then $\Mb$ is
free as an $R$--module with basis
$\set{m_{\s\t}|\s\in\Std_b(\lambda)\And\t\in\Std(\lambda)
                     \ForSome\lambda\in\Lab\cup\barLab}.$

\end{Lemma}

\proof As noted in (\ref{semistandard basis of Mmu}),
$\Mb$ is free as an $R$--module with basis $m_{\S\t}$,
where $\S\in\SStd(\lambda,\omega_b)$ and $\t\in\Std(\lambda)$ for
some multipartition $\lambda$ of $n$. It follows from
the definitions that $\SStd(\lambda,\omega_b)$ is non--empty if and
only if $\lambda\in\Lab\cup\barLab$; further, if 
$\lambda\in\Lab\cup\barLab$ then there is a bijection between
$\Std_b(\lambda)$ and $\SStd(\lambda,\omega_b)$ given by
$\s\mapsto\omega_b(\s)$ (\cf (\ref{type omega})). Consequently, if
$\S\in\SStd(\lambda,\omega_b)$ then $\S=\omega_b(\s)$ for a
uniquely determined $\s\in\Std_b(\lambda)$ and
$m_{\S\t}=m_{\s\t}$. Combining these observations gives the result.
\label{M^b basis}\endproof

Next we identify some elements in the kernel of $\theta_b$;
shortly we will see that these elements are actually a basis of
$\ker\theta_b$.

\begin{Lemma} Suppose that $0\le b\le n$ and $\lambda\in\barLab$. Let
$\s$ and $\t$ be standard $\lambda$--tableaux with
$\s\in\Std_b(\lambda)$. Then $\theta_b(m_{\s\t})=0$.
\label{kernel}

\end{Lemma}

\proof Suppose $h=\theta_b(m_{\s\t})\neq 0$. Let
$c=|\lambda^{(1)}|+\dots+|\lambda^{(s)}|$; then
$\lambda\in\Lambda_c^+$ \linebreak so that $b<c$. Now, because
$\lambda\in\Lambda_c^+$ and $\s\in\Std_b(\lambda)$, we can find a
permutation $w\in\Sym(\{b+1,\ldots,n\})$ such that $\s w\gedom\s$
and $\s w$ is a standard $\lambda$--tableau with the numbers
$1,\ldots,c$ all appearing in the first $s$ components of~$\s_w$.
(If $a_1<\dots<a_k$ are the numbers between $b+1$ and $c$ which
appear in the last $r-s$ components of~$\s$ and $b_1<\dots<b_k$ are
the $k$ smallest numbers larger than $b$ which appear in the first
$s$ components of $\s$ then we can set $w=(a_1,b_1)\dots(a_k,b_k)$.)
Let $\s_w=\s w$; then $d(\s)=d(\s_w)w^{-1}$ and
$m_{\s\t}=T_wm_{\s_w\t}$. Also let 
$\tilde w=w_{n-b,b}^{-1}w w_{n-b,b}$; then $\tilde w\in\Sym_{n-b}$ and
$\Twbn T_w=T_{\tilde w}\Twbn$ by Lemma~\ref{w_b properties 2}. Therefore, 
\begin{align*}
\theta_b(m_{\s\t})&=u_{n-b}^-T_{w_{n-b,b}}m_{\s\t}
        =u_{n-b}^-T_{w_{n-b,b}}T_wm_{\s_w\t}
        =u_{n-b}^- T_{\tilde w}\Twbn m_{\s_w\t}.\\
\intertext{Now, $\tilde w\in\Sym_{n-b}$ so 
$u_{n-b}^- T_{\tilde w}=T_{\tilde w}u_{n-b}^-$ by (\ref{L-comm})(iv);
consequently} 
\theta_b(m_{\s\t})&=T_{\tilde w}u_{n-b}^-T_{w_{n-b,b}} m_{\s_w\t}
       =T_{\tilde w}\theta_b(m_{\s_w\t}).
\end{align*}
Since $T_{\tilde{w}}$ is invertible, $\theta_b(m_{\s\t})\neq 0$ if
and only if $T_{\tilde{w}}\theta_b(m_{\s_w\t})\neq0$; therefore, we
may assume that $\s=\s_w$. Thus, it is enough to show that
$\theta(m_{\s\t})=0$ whenever $\s\in\Std_{c,n-c}(\lambda)$ and
$\t\in\Std(\lambda)$. However, if
$\s\in\Std_{c,n-c}(\lambda)$ then $m_{\s\t}\in M^{\omega_c}$  by
Lemma~\ref{M^b basis}; so $m_{\s\t}=u_{\omega_c}^+h$ for some $h\in\H$.
Now $u^+_{\omega_c}=u_c^+$;
so~$ \theta_b(m_{\s\t}) = u_{n-b}^-\Twbn m_{\s\t} 
        = u_{n-b}^-\Twbn u_c^+ h=0,
$ 
by Corollary~\ref{ubTucnull}, as desired.   
\endproof

Recall from Theorem~\ref{d=2} that given an integer $s$, with 
$1\le s\le r$, we let
$$f_s(q,\Q)=\prod_{1\le i\le s<j\le r}\prod_{-n<a<n}(q^aQ_i-Q_ja).$$
{\it For the remainder of this paper we assume that $f_s(q,\Q)$ is an 
invertible element of~$R$.}

Suppose that $\s\in\Std_{b,n-b}(\lambda)$ for some multipartition
$\lambda$. Then $\s'=\s w_{b,n-b}$ is a standard
$\lambda$--tableau which has the numbers $n-b+1,\dots,n$
appearing in its first~$s$ components and the remaining numbers
$1,\dots,n-b$ appearing in its last $r-s$ components.

\begin{Lemma} Suppose that $\lambda\in\Lab$ and let $\s$ and $\t$ be
standard $\lambda$--tableaux with $\s\in\Std_{b,n-b}(\lambda)$.
Let $\s'=\s w_{b,n-b}$. Then 
\begin{enumerate}\item $\Twbn m_{\s\t}=m_{\s'\t}$; and,
\item there exists an invertible element $\alpha\in R$ such that
$$\theta_b(m_{\s\t})\equiv \alpha m_{\s'\t} 
              +\sum_{\u'\gdom\s'} a_{\u'}m_{\u'\t}\mod\Nlambar,$$ 
for some $a_{\u'}\in R$. 
\end{enumerate}\label{independence}

\end{Lemma}

\proof (i) First note that $d(\s)\in\Sym_{(b,n-b)}$ and recall that
$w_{b,n-b}$ is a distinguished right coset representative of
$\Sym_{(b,n-b)}$. Therefore, 
$\len(d(\s'))=\len(d(\s))+\len(w_{b,n-b})$; consequently, 
$\Twbn m_{\s\t}=m_{\s'\t}$ since $w_{n-b,b}^{-1}=w_{b,n-b}$. 

(ii) Part (i) together with (\ref{L_k action}) implies that,
modulo $\Nlambar$,
\begin{align*}
\theta_b(m_{\s\t})&=u_{n-b}^-\Twbn m_{\s\t} 
   =\Big(\Prod_{t=1}^s(L_1-Q_t)\dots(L_{n-b}-Q_t)\Big)m_{\s'\t}\\
   &\equiv\Big(\Prod_{t=1}^s(\res_{\s'}(1)-Q_t)\dots
            (\res_{\s'}(n-b)-Q_t)\Big)m_{\s'\t}
          +\sum_{\u'\gdom\s'}a_{\u'}m_{\u'\t}
\end{align*}
for some $a_{\u'}\in R$. Let
$\alpha=\prod_{t=1}^s(\res_{\s'}(1)-Q_t)\dots(\res_{\s'}(n-b)-Q_t)$;
then $\alpha$ is the coefficient of $m_{\s'\t}$ in
$\theta_b(m_{\s\t})$. Now, $1,2,\dots,n-b$ all belong to one of
the last $r-s$ components of $\s'$; so, for $1\le k\le n-b$,
$\res_{\s'}(k)=q^jQ_c$ for some~$c$ and~$j$ with $s< c\le r$ and
$-n<j<n$.  Therefore, $\alpha$ is a product of terms of the form
$(q^jQ_c-Q_t)$, with $1\le t\le s$. As each of these factors
divides~$f_s(q,\Q)$, it follows that $\alpha$ is invertible; so
the Lemma is proved.
\endproof

Observe that the invertible element $\alpha$ in part (ii) of the
Lemma depends only on $\lambda=\Shape(\s)$, rather than $\s$ itself.
Notice also that 
$$ v_b=u_{n-b}^-\Twbn m_{\omega_b}
      =u_{n-b}^-\Twbn m_{\t^{\omega_b}\t^{\omega_b}}
      \equiv\alpha m_{\t\t^{\omega_b}}\mod\Nlambar,
$$ 
where $\t=\t^{\omega_b}w_{b,n-b}$; in particular, we have finally
proved that $v_b$ is non--zero. 

\begin{Definition} Suppose that $0\le b\le n$ and let
$\s\in\Std_{b,n-b}(\lambda)$ and $\t\in\Std(\lambda)$ for some
multipartition $\lambda\in\Lab$. Let
$v_{\s\t}=\theta_b(m_{\s\t})=u_{n-b}^-\Twbn m_{\s\t}$.

\end{Definition}

It follows from part (ii) of the Lemma that the $v_{\s\t}$ are
linearly independent elements in~$V^b$. In fact, they are a basis
of $V^b$.

\begin{Theorem} Suppose that $f_s(q,\Q)$ is invertible in
$R$ and let $b$ be an integer with $0\le b\le n$. Then 
$V^b$ is free as an $R$--module with basis
$$\set{v_{\s\t}|\s\in\Std_{b,n-b}(\lambda)\And\t\in\Std(\lambda)
                 \ForSome\lambda\in\Lab}$$
and $\ker\theta_b$ is free as an $R$--module with basis
$$\set{m_{\u\v}|\u\in\Std_{b}(\mu)\And\v\in\Std(\mu)
                 \ForSome\mu\in\barLab}.$$ 
\label{V^b basis}

\end{Theorem}

\proof The homomorphism $\theta_b$ is surjective, so by 
Lemma~\ref{M^b basis} $V^b$ is spanned by the
elements $\theta_b(m_{\u\v})$, where $\u\in\Std_b(\mu)$,
$\v\in\Std(\mu)$ and $\mu\in\Lab\cup\barLab$. Furthermore, 
$\theta_b(m_{\u\v})=0$ whenever $\mu\in\barLab$, by
Lemma~\ref{kernel}. Finally, the elements 
$\set{v_{\s\t}=\theta_b(m_{\s\t})|\s\in\Std_{b,n-b}(\lambda),
                     \t\in\Std(\lambda)\And\lambda\in\Lab}$
are linearly independent (and non--zero) by Lemma~\ref{independence}(ii).
Combing these three statements proves the Theorem.
\endproof

Let $\Nb=\bigcap_{\lambda\in\Lab}\Nlambar
               =\sum_{\mu\in\barLab}N^\mu$. 
Then $\Nb$ is a two--sided ideal in $\H$ and it is free 
as an $R$--module with basis
$\set{m_{\u\v}|\u,\v\in\Std(\mu)\ForSome\mu\in\barLab}$.

\begin{Corollary} Suppose that $0\le b\le n$. Then
$\ker\theta_b=\Mb\cap\Nb$.
\label{N_b quotient}
\end{Corollary}

\begin{Remark} The proof of Theorem~\ref{V^b basis} relies on the
assumption that $f_s(q,\Q)$ is an invertible element of $R$.
Without assuming that $f_s(q,\Q)$ is invertible it is possible to
prove via a brute force calculation that~$V^b$ is free as an $R$--module 
with basis
$$\Set[55]v_bL_1^{d_1}\dots L_n^{d_n}T_w|$w\in\Sym_n$,
         $0\le d_i<s$ for $1\le i\le b$\\ and
         $0\le d_i<r-s$ for $b<i\le n$|.$$
(A straightforward argument using the Robinson--Schensted
correspondence verifies that the rank of $V^b$ is the same in both
cases.) That this set is a basis of~$V^b$ can also be deduced from
Proposition~\ref{proof of gensit} below; however, this argument requires
that $f_s(q,\Q)$ be invertible in $R$.
\label{gensit}
\end{Remark}

In \cite[Corollary~4.15]{DJM:cyc} Specht filtrations of the right
ideals $M^\mu$ of $\H$ were constructed as a
consequence of (\ref{semistandard basis of Mmu}); we now refine the
filtration of $\Mb$ to give Specht filtrations of the modules
$V^b$ and $\ker\theta_b$.

\begin{Theorem} Suppose that $f_s(q,\Q)$ is invertible in
$R$ and let $b$ be an integer with $0\le b\le n$.
\begin{enumerate}\item There is a filtration 
$V^b=V_1\supset V_2\supset\dots\supset V_k\supset V_{k+1}=0$ of $V^b$
such that for each $1\le i\le k$ there exists a multipartition
$\lambda_i\in\Lab$ with $V_i/V_{i+1}\cong S^{\lambda_i}$.
Moreover, for each $\lambda\in\Lab$ the number of $i$ with
$\lambda_i=\lambda$ is $|\Std_{b,n-b}(\lambda)|$.
\item There is a filtration 
$\ker\theta_b=K_1\supset K_2\supset\dots\supset K_l\supset K_{l+1}=0$ of $\ker\theta_b$
such that for each $1\le i\le l$ there exists a multipartition
$\mu_i\in\barLab$ with $K_i/K_{i+1}\cong S^{\mu_i}$.
Moreover, for each $\mu\in\barLab$ the number of $i$ with
$\mu_i=\mu$ is $|\Std_{b}(\mu)|$.
\end{enumerate}\label{Specht filtration}

\end{Theorem}

\proof We recall the construction from \cite[Cor.~4.15]{DJM:cyc},
using Lemma~\ref{M^b basis} to adapt the notation. Let 
$\s_1,\s_2,\dots,\s_N$ be the tableaux in
$\bigcup_{\lambda\in\Lab\cup\barLab}\Std_b(\lambda)$, ordered so
that $j>i$ whenever $\lambda_i\gdom\lambda_j$; here we set
$\lambda_i=\Shape(\s_i)$ for all $i$. For $i=1,2,\dots,N$ let
$M_i$ be the $R$--submodule of $\Mb$ with basis
$\set{m_{\s_j\t}|j\ge i\And\t\in\Std(\lambda_j)}.$
Then the proof of \cite[Cor.~4.15]{DJM:cyc} shows that each $M_i$ is
an $\H$--submodule of $\Mb$ and that $M_i/M_{i+1}\cong
S^{\lambda_i}$ for all $i$.

Choose $k$ to be maximal such that $\lambda_k\in\Lab$; then
$M_{k+1}=\ker\theta_b$ and $\lambda_i\in\barLab$ if and only if
$i>k$. By Corollary~\ref{N_b quotient} we can set $V_i=M_i/\ker\theta_b$,
for $i=1,2,\dots,k+1$, to obtain a filtration of $V^b$ with the
required properties. Similarly, setting $K_j=M_{k+j}$, for
$j=1,2,\dots,N-k$, gives the promised filtration of
$\ker\theta_b$.
\endproof

Recall from before Corollary~\ref{content} that the content $\cont(\lambda)$
of $\lambda$ is the sequence $(c_r)_{r\in R}$, where
$c_r=\#\set{x\in[\lambda]|\res(x)=r}$ for all $r\in R$.
Importantly, two Specht modules belong to the same block only if the
corresponding multipartitions have the same content by Corollary~\ref{content}.

\begin{Corollary} Suppose that $f_s(q,\Q)$ is invertible and let $b$ and $c$
be distinct integers with $0\le b,c\le n$. Then
$\Hom_\H(V^b,V^c)=0$.
\label{vanishing homs}

\end{Corollary}

\proof Let $\mathcal R_s=\set{q^dQ_i|-n<d<n\And 1\le i\le s}$; then
$\mathcal R_s$ is the complete set of possible residues $\res(x)$, where
$x$ runs through the nodes in the first $s$--components of the diagram of
any multipartition of $n$.  In addition, since $f_s(q,\Q)$ is invertible, if
$y$ is a node appearing in one of the last $r-s$ components of some
multipartition then $\res(y)\notin\mathcal R_s$; consequently,
if~$\lambda$ is a multipartition and $\cont(\lambda)=(c_r)_{r\in R}$ then 
$\lambda\in\Lab$ if and only if $b=\sum_{r\in\mathcal R_s}c_r$. Therefore, if
$\lambda\in\Lab$ and $\mu\in\Lambda_c^+$ then $\cont(\lambda)\ne\cont(\mu)$
--- note that by assumption $b\ne c$.

Now consider $\Hom_\H(V^b,V^c)$. By Theorem~\ref{Specht filtration}, $V^b$
has a Specht filtration indexed by multipartitions in $\Lab$ and
$V^c$ has a Specht filtration indexed by the multipartitions in
$\Lambda_c^+$. Therefore, the simple composition factors of~$V^b$
and~$V^c$ belong to different blocks by the last paragraph and
Corollary~\ref{content}; hence, $\Hom_\H(V^b,V^c)=0$ by Schur's lemma.
\endproof

Notice that $\barLab=\bigcup_{c=b+1}^n\Lambda_c^+$; therefore, by
a similar argument, again using Theorem~\ref{Specht filtration} and
Corollary~\ref{content}, the composition factors of $V^b$ and
$\ker\theta_b$ belong to different blocks.  Hence, we also have
the following.

\begin{Corollary} Suppose $f_s(q,\Q)$ is invertible in $R$
and that $0\le b\le n$. Then the composition factors of~$V^b$ and 
$\ker\theta_b$ belong to different blocks of $\H$; consequently,
$\Mb\cong V^b\oplus\ker\theta_b$ and
$\End_\H(\Mb)\cong \End_\H(V^b)\oplus\End_\H(\ker\theta_b).$
\label{splitting}
\end{Corollary}

In fact, this allows us to determine a basis of $\End_\H(V^b)$.
Given two standard $\lambda$--tableaux $\s$ and $\t$ in
$\Std_{b,n-b}(\lambda)$ let $\theta_{\s\t}\map{V^b}V^b$ be the
$R$--linear map given by $\theta_{\s\t}(v_bh)=v_{\s\t}h$ for all
$h\in\H$; then $\theta_{\s\t}$ is an $R$--module homomorphism. 
{\it A priori} there is no reason to expect that $\theta_{\s\t}$ is even
well defined; nevertheless, it is and these elements give a basis of
$\End_\H(V^b)$.

\begin{Theorem}  Suppose $f_s(q,\Q)$ is invertible in $R$
and let $0\le b\le n$. Then $\End_\H(V^b)$ is free as an
$R$--module with basis 
$$\set{\theta_{\s\t}|\s,\t\in\Std_{b,n-b}(\lambda)
                    \ForSome\lambda\in\Lab}.$$
\label{hom-space basis}

\end{Theorem}

\proof A basis of $\End_\H(\Mb)$ is given by (\ref{SStd basis});
in light of Lemma~\ref{M^b basis}, we see that $\End_\H(\Mb)$ has as
basis the maps 
$$\set{\phi_{\s\t}|\s,\t\in\Std_b(\lambda)
                    \ForSome\lambda\in\Lab\cup\barLab},$$
where $\phi_{\s\t}$ is given by
$\phi_{\s\t}(m_{\omega_b}h)=m_{\s\t}h$ for all $h\in\H$.

By Corollary~\ref{splitting} the homomorphism $\theta_b\map\Mb V^b$ splits, so
$\theta_b$ has a right inverse; we abuse notation and write $\theta_b^{-1}$
for this one sided inverse. Let~$\phi$ be any map in
$\End_\H(\Mb)$; then $\theta_b\phi\theta_b^{-1}$ belongs to $\End_\H(V^b)$
and every homomorphism in $\End_\H(V^b)$ is of this form. Now,
$\theta_b(m_{\omega_b})=v_b$; so there exists an $h_b\in\ker\theta_b$ such
that $\theta_b^{-1}(v_bh)=(m_{\omega_b}+h_b)h$ for all $h\in\H$.  Observe
that $\phi_{\s\t}(h_b)\in\ker\theta_b$ since $\ker\theta_b$ are $V^b$ are
in different blocks by Corollary~\ref{splitting}. Therefore, by Theorem~\ref{V^b basis},
$$\theta_b\phi_{\s\t}\theta_b^{-1}(v_bh)
      =\theta_b\phi_{\s\t}(m_{\omega_b}+h_b)h
      =\theta_b(m_{\s\t})h
      =\cases v_{\s\t}h,&\!\!\If\lambda\in\Lab,\\
              0,&\!\!\If \lambda\in\barLab.
       \endcases$$
Consequently, $\theta_b\phi_{\s\t}\theta_b^{-1}=\theta_{\s\t}$ if
$\lambda\in\Lab$ and $\theta_b\phi_{\s\t}\theta_b^{-1}=0$ if
$\lambda\in\barLab$; in particular, $\theta_{\s\t}\in\End_\H(V^b)$ whenever
$\s$ and $\t$ belong to $\Std_{b,n-b}(\lambda)$. As the elements
$\set{v_{\s\t}|\s,\t\in\Std_{b,n-b}(\lambda)\ForSome\lambda\in\Lab}$ are
linearly independent, so are the corresponding homomorphisms
$\{\theta_{\s\t}\}$. The Theorem follows.
\endproof

The argument also shows that a basis of $\End_\H(\ker\theta_b)$
is given by restricting the maps 
$\set{\phi_{\s\t}|\s,\t\in\Std_b(\mu)\ForSome\mu\in\barLab}$ to 
$\ker\theta_b$.

The final property of the modules $V^b$ that we need is that they
are projective $\H$--modules. To prove this let $\D_{(b,n-b)}$ be
the set of distinguished right coset representatives of
$\Sym_{(b,n-b)}$ in $\Sym_n$; see, for
example,~\cite[Prop.~3.3]{M:ULect}. 

\begin{Theorem} Suppose that $f_s(q,\Q)$ is invertible. Then
$\DS\H\cong\bigoplus_{b=0}^n \tbinom n b V^b$.
\label{projective decomposition}

\end{Theorem}

\proof Suppose that $0\le b\le n$. As in the proof of 
Theorem~\ref{hom-space basis} let $\theta_b^{-1}$ be a right inverse to
$\theta_b$ and set $V_b'=\theta_b^{-1}(V^b)$; then $V_b'\cong V^b$
is a submodule of $\Mb$. For each multipartition $\lambda\in\Lab$
and tableaux $\s\in\Std_{b,n-b}(\lambda)$ and $\t\in\Std(\lambda)$
let $v_{\s\t}'=\theta_b^{-1}(v_{\s\t})$; then $v_{\s\t}'\in V_b'$
and $v_{\s\t}'=m_{\s\t}+h_{\s\t}$ for some
$h_{\s\t}\in\ker\theta_b$.  Moreover, the set of these elements is a
basis of~$V_b'$.

We claim that
$\H=\bigoplus_{b=0}^n\bigoplus_{w\in\D_{(n-b,b)}}  T_w^*V_b';$ 
since $[\Sym_n:\Sym_{(b,n-b)}]=\binom n b$ this will establish the
Theorem. Note that $T_w^*V_b'\cong V^b$ as right $\H$--modules.

If $w\in\D_{(b,n-b)}$ and $\s\in\Std_{b,n-b}(\lambda)$ then $\s_w=\s w$
is again a standard $\lambda$--tableau, $d(\s_w)=d(\s)w$ and
$\len(d(\s_w))=\len(d(\s))+\len(w)$.  Therefore, $T_w^*V_b'$ has
as basis the elements
$$T_w^*v_{\s\t}'=T_w^*\(m_{\s\t}+h_{\s\t}\)
                =m_{\s_w\t}+T_w^*h_{\s\t},$$
where $\s,\t\in\Std_{b,n-b}(\lambda)$ for some $\lambda\in\Lab$.
Notice that $\ker\theta_b\subseteq\Nb$ by Corollary~\ref{N_b quotient} and that 
$\Nb$ is a two--sided ideal of $\H$; so $T_w^*h_{\s\t}\in\Nb$. 

Now, if $\lambda\in\Lab$ then
$\Std(\lambda)=\coprod_{w\in\D_{(b,n-b)}}\Std_{b,n-b}(\lambda)w$.
Therefore, the elements 
$$\bigcup_{b=0}^n\set{T_w^*v_{\s\t}'|w\in\D_{(b,n-b)},
        \s\in\Std_{b,n-b}(\lambda)
        \And\t\in\Std(\lambda)\ForSome\lambda\in\Lab}$$ 
are linearly independent; moreover, by (\ref{standard basis}) and the
above remarks, this set is also a basis of $\H$. Consequently,
$\H=\bigoplus_{b=0}^n\bigoplus_{w\in\D_{(b,n-b)}} T_w^*V_b'$, proving
our claim and hence the Theorem.
\endproof

A similar argument shows that 
$M^{\omega_a}\cong\bigoplus_{b=a}^n\tbinom{n-a}bV^b$ whenever 
$0\le a\le n$; the Theorem is the case $a=0$ since $m_{\omega_0}=1$.

Recall that a {\sf progenerator}, or projective generator, for $\H$
is a projective $\H$--module which contains every principal
indecomposable $\H$--module as a direct summand.  Every direct
summand of $\H$ is automatically projective and every principal
indecomposable $\H$--module is a direct summand of $\H$;
hence, Theorem~\ref{projective decomposition} yields the following.

\begin{Corollary} Suppose that $f_s(q,\Q)$ is invertible. Then
\begin{enumerate}\item$V^b$ is a projective $\H$--module for $b=0,1,\dots,n$; and,
\item $\DS V=\bigoplus_{b=0}^n V^b$ is a progenerator for $\H$.
\end{enumerate}\label{progenerator}
\end{Corollary}

\section{The Morita equivalence}

We are almost ready to prove our main results; in fact, we have
already done most of the hard work. In Corollary~\ref{progenerator}(ii) we
showed that $V=\bigoplus_{b=0}^n V^b$ is a progenerator of $\H$;
hence, it follows that $\H$ is Morita equivalent to $\End_\H(V)$.
To prove Theorem~\ref{Morita equivalence} we show that $\End_\H(V)$ is
isomorphic to a direct sum of tensor products of smaller
Ariki--Koike algebras.

\begin{Proposition} Suppose that $f_s(q,\Q)$ is invertible. Then 
$\H$ is Morita equivalent to $\bigoplus_{b=0}^n\End_\H(V^b)$.
\label{morita}

\end{Proposition}

\proof Let $V=\bigoplus_{b=0}^n V^b$. As remarked above, it follows from
Corollary~\ref{progenerator}(ii) that $\H$ is Morita equivalent to $\End_\H(V)$; see
\cite[Lemma~2.2.3]{Benson:I}. However, by Corollary~\ref{vanishing homs}, 
$$\End_\H(V)=\bigoplus_{0\le b,c\le n}\Hom_\H(V^b,V^c)
            \cong\bigoplus_{b=0}^n \End_\H(V^b),$$
giving the result.
\endproof

For $b=0,1\dots,n$ let $\Hb=\H_{q,\Q_1}(b)\otimes\H_{q,\Q_2}(n-b)$,
where $\Q_1=(Q_1,\dots,Q_s)$ and $\Q_2=(Q_{s+1},\dots,Q_r)$. We have
already computed $\End_\H(V^b)$ in Theorem~\ref{hom-space basis}; in order
to prove Theorem~\ref{d=2}, and hence our main result, we will use this
result to show that $\End_\H(V^b)\cong\Hb$.

The subalgebra $\H(\Sym_b)\otimes\H(\Sym_{n-b})$ of $\Hb$ is
isomorphic to the subalgebra $\H(\Sym_{(n-b,b)})$ of $\H$ spanned by
$\set{T_w|w\in\Sym_{(n-b,b)}}$ --- note that $\H(\Sym_{(b,n-b)})$
and $\H(\Sym_{(n-b,b)})$ are isomorphic algebras, the reason for
introducing the twist is that
$\H(\Sym_{(n-b,b)})v_b=v_b\H(\Sym_{(b,n-b)})$ by 
Proposition~\ref{v_b properties}. In general $\H$ has no subalgebra isomorphic
to $\H_b\otimes\H_{n-b}$; however, $\H$ does have an $R$--submodule
isomorphic to $\Hb$ and this we can exploit.  Let $\H_{n-b,b}$ be
the $R$--submodule of $\H$ spanned by the elements
$$\Set[72]L_1^{d_1}\dots L_n^{d_n}T_w| $w\in\Sym_{(n-b,b)}$,\ 
         $0\le d_i<r-s$ for $1\le i\le n-b$\\ and
         $0\le d_i<s$ for $n-b<i\le n$|.$$
We emphasize that typically $\H_{n-b,b}$ is not a subalgebra of $\H$;
however, $\H_{n-b,b}$ does generate the subalgebra $H_b$ of Corollary~\ref{endhom}.

\begin{Point}* Suppose that $0\le b\le n$. Then $\Hb$ and $\H_{n-b,b}$
are isomorphic as $R$--modules via the $R$--linear map
$\Theta_b\map\Hb\H_{n-b,b}$ determined by
$$L_1^{d_1}\dots L_b^{d_b}T_x\otimes 
  L_1^{e_1} \dots L_{n-b}^{e_{n-b}}T_y 
 \mapsto (L_1^{e_1}\dots L_{n-b}^{e_{n-b}}T_y)
         (L_{n-b+1}^{d_1}\dots L_n^{d_b}T_{x'}),$$
where $x'=w_{n-b,b}xw_{b,n-b}=w^{-1}_{b,n-b}xw_{b,n-b}$;
here $x\in\Sym_b$, $y\in\Sym_{n-b}$, $0\le d_i<s$ for $i=1,\dots,b$
and $0\le e_j<r-s$ for $j=1,\dots,n-b$.
\label{identification}
\end{Point}

\noindent Maintaining the notation of (\ref{identification}) notice that 
$$(L_1^{e_1}\dots L_{n-b}^{e_{n-b}}T_y)
         (L_{n-b+1}^{d_1}\dots L_n^{d_b}T_{x'})
   =L_1^{e_1}\dots L_{n-b}^{e_{n-b}}L_{n-b+1}^{d_1}\dots L_n^{d_b}
      T_{x'y}$$
by (\ref{L-comm}); in particular, the map
$\Sym_b\times\Sym_{n-b}\longrightarrow\Sym_{(n-b,b)}\,{:}\,(x,y)
        \mapsto x'y=yx'$
is an isomorphism of groups. Hereafter, we identify elements of
the $R$--modules $\Hb$ and $\H_{n-b,b}$ via the map $\Theta_b$.

\begin{Lemma} Suppose that $0\le b\le n$. Then $V^b$ becomes a left
$\Hb$--module via $hv=\Theta_b(h)v$ for all $h\in\Hb$ and 
$v\in V^b$. Consequently, 
$\Theta_b(h_1h_2)v=\Theta_b(h_1)\Theta_b(h_2)v$ for all 
$h_1,h_2\in\Hb$ and all $v\in V^b$.
\label{left H_b-module}

\end{Lemma}

\proof It is enough to check that the relations in $\Hb$ are
preserved.  As an algebra $\Hb$ is generated by the two sets of
commuting elements $T_i\otimes1$ and $1\otimes T_j$, where $0\le i<b$
and $0\le j<n-b$. Now
$$\Theta_b(T_i\otimes 1)=\cases L_{n-b+1},&\If i=0,\\
                                T_{n-b+i},&\Otherwise,
                         \endcases
\quad\And\quad
\Theta_b(1\otimes T_j)=\cases L_1,&\If j=0,\\
                              T_j,&\Otherwise.
                         \endcases
$$
Therefore, the relations in $\H$ and (\ref{L-comm})(ii) ensure that
all of the relations in $\Hb$ are satisfied except possibly that
$(T_0\otimes1-Q_1)\dots(T_0\otimes1-Q_s)$ and 
$(1\otimes T_0-Q_{s+1})\dots(1\otimes T_0-Q_r)$ must both act as
zero on~$V^b$. However, this is precisely the content of parts (ii) 
and (i), respectively, of Corollary~\ref{v_b annihilators}. (Notice that when
$\Theta_b$ is applied to the braid
relation $T_0T_1T_0T_1=T_1T_0T_1T_0$ in~$\H_b$ it becomes
$qL_{n-b+1}L_{n-b+2}=qL_{n-b+2}L_{n-b+1}$ in $\H$; this holds by
virtue of~(\ref{L-comm})(ii).)
\endproof

Therefore, $V^b$ is an $(\Hb,\H)$--bimodule. Evidently, the left action
of $\Hb$ on $V^b$ commutes with the right action of $\H$ so we have a
homomorphism from $\Hb$ into $\End_\H(V^b)$. Furthermore, 
Lemma~\ref{left H_b-module} implies that for all $h_1,h_2\in\Hb$ there exists
some $h_3\in\H$ such that 
$\Theta_b(h_1h_2)=\Theta_b(h_1)\Theta_b(h_2)+h_3$
and $h_3\in\mathsf L_{H_b}(V^b)$; \cf Corollary~\ref{endhom}.

Shortly we will see that $\Hb$ acts faithfully on $V^b$; before we
can show this we need some more notation. Recall that
$\Lambda^+=\Lamp(n:r)$ and, more generally, that $\Lamp(b:s)$ is the set of
multipartitions of $b$ with $s$--components. There is an evident bijection
$\Lamp(b:s)\times\Lamp(n-b:r-s)\longrightarrow\Lab$ which maps the pair
$\sigma=(\sigma^{(1)},\dots,\sigma^{(s)})\in\Lamp(b:s)$ and
$\tau=(\tau^{(1)},\dots,\tau^{(r-s)})\in\Lamp(n-b:r-s)$ to the
multipartition $(\sigma,\tau)=(\sigma^{(1)},\dots,\sigma^{(s)},
\tau^{(1)},\dots,\tau^{(r-s)})\in\Lab$.

\begin{Lemma} Suppose that $0\le b\le n$ and that $\lambda\in\Lab$.
Then $\lambda=(\sigma,\tau)$ for unique multipartitions
$\sigma\in\Lamp(b:s)$ and $\tau\in\Lamp(n-b:r-s)$; moreover,
we have $\theta_b(u_\lambda^+)=\Theta_b(u_\sigma^+\otimes u_\tau^+)v_b$.
\label{ulambdaactsright}

\end{Lemma}

\proof The uniqueness of $\sigma$ and $\tau$ follows from the
remarks above. For the remaining statement first recall that
$u_{\alpha,t}=\prod_{k=1}^\alpha(L_k-Q_t)$ for any $\alpha$ and $t$. Let
$\alpha_t=|\lambda^{(1)}|+\dots+|\lambda^{(t-1)}|$, for $1\le t\le r$,
and let $\beta_t=\alpha_{t+r-s}-b$ for $t=1,\dots,r-s$.
Then $u_\lambda^+=u_{\alpha_1,1}\dots u_{\alpha_r,r}$ and, abusing 
notation,
$u_\sigma^+\otimes u_\tau^+=u_{\alpha_1,1}\dots u_{\alpha_s,s}\otimes
             u_{\beta_1,s+1}\dots,u_{\beta_{r-s},r}\in\Hb$.

Now, by Proposition~\ref{v_b properties}(iv), if $1\le t\le s$ then
$$\Theta_b(u_{\alpha_t,t}\otimes 1)v_b
     =\Big(\prod_{k=n-b+1}^{n-b+\alpha_t}(L_k-Q_t)\Big)v_b
      =v_b u_{\alpha_t,t};$$
so, 
$\Theta_b(u_\sigma^+\otimes 1)v_b=v_b u_{\alpha_1,1}
          \dots u_{\alpha_s,s}$.
On the other hand, by Proposition~\ref{v_b properties}(iii),
$$\Theta_b(1\otimes u_{\beta_t,t})v_b
    =\Big(\prod_{k=1}^{\beta_t}(L_k-Q_t)\Big)v_b
     =v_b\prod_{k=b+1}^{b+\beta_t}(L_k-Q_t),$$
for $t=s+1,\dots,r$. Now,
$u_{\alpha_t,t}=\prod_{k=1}^{b+\beta_t}(L_k-Q_t)$ for $t=s+1,\dots,r$;
therefore,
\begin{align*}
\Theta_b(1\otimes u_\tau^+)v_b
    &=v_b \prod_{t=s+1}^r\prod_{k=b+1}^{b+\beta_t}(L_k-Q_t)\\
    &=u_{n-b}^-\Twbn u_b^+
      \prod_{t=s+1}^r\prod_{k=b+1}^{b+\beta_t}(L_k-Q_t)\\
    &=u_{n-b}^-\Twbn u_{\alpha_{s+1},s+1}\dots u_{\alpha_r,r}.
\end{align*}
Hence, $\Theta_b(u_\sigma^+\otimes u_\tau^+)v_b
      =u_{n-b}^-\Twbn u_{\alpha_1,1}\dots u_{\alpha_r,r}
      =\theta_b(u_\lambda^+)$
as claimed.
\endproof

We want to extend this result to the standard basis of $\Hb$.
Suppose that $\lambda\in\Lab$ and, as above, write
$\lambda=(\sigma,\tau)$ for multipartitions $\sigma\in\Lamp(b:s)$ and
$\tau\in\Lamp(n-b:r-s)$. Given standard tableaux
$\s_1=(\s_1^{(1)},\dots,\s_1^{(s)})\in\Std(\sigma)$
and
$\s_2=(\s_2^{(1)},\dots,\s_2^{(r-s)})\in\Std(\tau)$
let $(\s_1,\s_2')$ be the $\lambda$--tableau
$$(\s_1,\s_2')=(\s_1^{(1)},\dots,\s_1^{(s)},
       \s_2^{(1)}w_{n-b,b},\dots,\s_2^{(r-s)}w_{n-b,b});$$
then $(\s_1,\s_2')\in\Std_{b,n-b}(\lambda)$ and a straightforward
calculation reveals the following.
 
\begin{Lemma} Suppose that $\lambda=(\sigma,\tau)\in\Lab$ as above. Then 
the map
$$\Std(\sigma)\times\Std(\tau)\longrightarrow\Std_{b,n-b}(\lambda)
        \,{:}\,(\s_1,\s_2)\mapsto\s=(\s_1,\s_2')$$
is a bijection. Furthermore,
$d(\s)=d(\s_1)w_{n-b,b}^{-1}d(\s_2)w_{n-b,b}\in\Sym_{(b,n-b)}$
and $\len(d(\s))=\len(d(\s_1))+\len(d(\s_2))$.
\label{std decomp}
\end{Lemma}

We can now give the connection between $\theta_b$ and $\Theta_b$.

\begin{Lemma} Suppose $\lambda=(\sigma,\tau)\in\Lab$ as above. Let
$\s,\t\in\Std_{b,n-b}(\lambda)$ and write $\s=(\s_1,\s_2')$ and
$\t=(\t_1,\t_2')$ as in Lemma~\ref{std decomp}.  Then
$$\theta_b(m_{\s\t})=\Theta_b(m_{\s_1\t_1}\otimes m_{\s_2\t_2})v_b.$$
\label{nice action}

\end{Lemma}

\proof By Proposition~\ref{v_b properties},
$\H(\Sym_{(n-b,b)})v_b=v_b\H(\Sym_{(b,n-b)})$; therefore,
$$\Theta_b(x_\sigma\otimes x_\tau)v_b
      =x_{(\tau,\sigma)}v_b
      =v_bx_{(\sigma,\tau)}
      =v_bx_\lambda.$$
Furthermore, $\Theta_b(u_\sigma^+\otimes u_\tau^+)v_b=\theta_b(u_\lambda^+)$
by Lemma~\ref{ulambdaactsright}; therefore, using Lemma~\ref{left H_b-module} to combine 
these two equations shows that
\begin{align*}
\Theta_b(m_\sigma\otimes m_\tau)v_b
  &=\Theta_b(x_\sigma u_\sigma^+\otimes x_\tau u_\tau^+)v_b
   =\Theta_b(x_\sigma\otimes x_\tau)\Theta_b(u_\sigma^+\otimes u_\tau^+)v_b\\
  &=\Theta_b(x_\sigma\otimes x_\tau)v_bu_\lambda^+
   =v_b x_\lambda u_\lambda^+
   =v_b m_\lambda.
\end{align*}
Finally, note that if $x\in\Sym_b$ and $y\in\Sym_{n-b}$ then, by
(\ref{identification}) and Proposition~\ref{v_b properties},
$$\Theta_b(T_x\otimes T_y)v_b
       =T_{w_{b,n-b}^{-1}xw_{b,n-b}}T_yv_b
       =v_bT_xT_{w_{n-b,b}^{-1}yw_{n-b,b}}.$$
Since
$m_{\s_1\t_1}\otimes m_{\s_2\t_2}
     =(T_{d(\s_1)}^*\otimes T_{d(\t_1)}^{\phantom*})(m_\sigma\otimes m_\tau)
         (T_{d(\s_2)}^*\otimes T_{d(\t_2)}^{\phantom*}),$
another application of Lemma~\ref{left H_b-module}, together with
Lemma~\ref{std decomp}, now completes the proof.
\endproof
 
We have finally reached the summit.

\begin{Theorem} Suppose that $f_s(q,\Q)$ is invertible in $R$
and let $b$ be an integer with $0\le b\le n$. Then
\begin{enumerate}\item $V^b$ is a faithful left $\H_b\otimes\H_{n-b}$--module;
\item $\End_\H(V^b)\iso\Hb$; and,
\item $\H$ is Morita equivalent to $\DS\H_s=\bigoplus_{b=0}^n\Hb$.
\end{enumerate}\label{endomorphism ring}

\end{Theorem}

\proof By Lemma~\ref{left H_b-module} $\Hb$ acts on $V^b$ by left multiplication;
this action is faithful because, by Lemma~\ref{nice action} and Theorem~\ref{V^b basis},
the elements of the standard basis of $\Hb$ map $v_b$ to linearly
independent elements in $V^b$. Further, because the action is faithful,
$\End_\H(V^b)\iso\Hb$ by Theorem~\ref{hom-space basis} and Lemma~\ref{nice action}.
Finally, part~(iii) follows from part~(ii) and Proposition~\ref{morita}.
\endproof

This proves Theorem~\ref{d=2} and hence the main result of this paper. The
remainder of this section examines this Morita equivalence 
more closely.

For an algebra $A$ let $\Mod_A$ be the category of (finite dimensional)
right $A$--modules. Then we have shown that there exist functors (in fact,
category equivalences),
$$\Ho_s\map{\Mod_\H}\Mod_{\H_s}\And
  \hatHo_s\map{\Mod_{\H_s}}\Mod_\H\!.$$
These functors are described in terms of the 
$(\H_s,\H)$--bimodule $V=\bigoplus_{b=0}^nV^b$; explicitly,
$\Ho_s(M)=\Hom_{\H}(V,M)$ and $\hatHo_s(X)=X\otimes_{\H_s}V$.
(Here we consider $V$ as a left $\H_s$--module by specifying that
$\Hb$ annihilates $V^c$ when $b\ne c$.) For this, and other standard
facts about Morita equivalences, see \cite[Section~2.2]{Benson:I}.

By Theorem~\ref{projective decomposition} we can write
$\H=\bigoplus_{b=0}^n\H(b)$ where $\H(b)$ is the smallest two--sided
ideal of $\H$ which contains $V^b$ as a direct summand (for 
$0\le b\le n$).  Therefore, $\Mod_\H=\bigoplus_{b=0}^n\Modb$. The
modules in $\Modb$ are mapped into $\Mod_{\Hb}$ by $\Ho_s$, so we
have subfunctors
$$\Ho_{s,b}\map{\Modb}\Mod_{\Hb}\And
  \hatHo_{s,b}\map{\Mod_{\Hb}}\Modb$$
given by $\Ho_{s,b}(M)=\Hom_{\H}(V^b,M)$ and
$\hatHo_{s,b}(X)=X\otimes_{\Hb}V^b$. Each of these functors induces a
Morita equivalence. 

In general, $V$ is not free as a left $\H_s$--module; however,  as
the next result shows,  $V^b$ is free as a left $\Hb$--module.
Recall that $\D_{(b,n-b)}$ is the set of distinguished right coset
representatives of $\Sym_{(b,n-b)}\cong \Sym_b\times\Sym_{n-b}$ in
$\Sym_n$. For each $w\in\D_{(b,n-b)}$ let $V^b_w$ be the
$R$--submodule $\Theta_b(\Hb)v_bT_w$ of $V^b$.

\begin{Proposition} Suppose that $0\leq b\leq n$ and let $w\in\D_{(b,n-b)}$. Then
$V^b_d$ is a left $\Hb$--submodule of $V^b$ which is isomorphic to
the left regular representation of $\Hb$. A basis of $V^b_d$ is
given by 
$$\set{v_{\s\t}|\s\in\Std_{b,n-b}(\lambda)\And\t\in\Std(\lambda)w
                 \ForSome\lambda\in\Lab}.$$
Moreover, as a left $\Hb$--module,
$$\DS V^b=\bigoplus_{w\in\D_{(b,n-b)}}V^b_w.$$
Consequently, $V^b$ is a free $\Hb$--module of rank $\binom n b$.
\label{proof of gensit}

\end{Proposition}

\proof Applying Lemma~\ref{nice action} and Theorem~\ref{V^b basis} we see that
$V^b_w$ has basis
$\set{v_{\s\t}T_w|\s,\t\in\Std_{b,n-b}(\lambda)\ForSome\lambda\in\Lab}$.
Now, if $\t\in\Std_{b,n-b}(\lambda)$ then $d(\t)\in\Sym_{b,n-b}$ so
that $\len(d(\t)w)=\len(d(\t))=\len(w)$ since $w$ is a distinguished
right coset representative of $\Sym_{(b,n-b)}$ in $\Sym_n$.
Therefore, if $\s$ and $\t$ in
$\Std_{b,n-b}(\lambda)$ then $m_{\s\t}T_w=m_{\s\u}$ where 
$\u=\t w\in\Std(\lambda)$; so, $v_{\s\t}T_w=v_{\s\u}$. Therefore,
$V^b_w$ has the required basis and
$V^b=\bigoplus_{w\in\D_{(b,n-b)}}V^b_w$.  Finally, $V^b_w$ affords
the left regular representation of $\Hb$ by 
Theorem~\ref{endomorphism ring}(i).
\endproof

Suppose that $\set{h_i| i\in I}$ is a basis of $\Hb$; then, by the
Proposition, 
$$\set{\Theta_b(h_i)v_bT_w|i\in I\And w\in\D_{(b,n-b)}}$$
is a basis of $V^b$. In particular, taking $\set{h_i|i\in I}$ to be the
Ariki--Koike basis of $\Hb$ we obtain the basis of $V^b$ mentioned in
Remark~\ref{gensit}.

In particular, since  $V^b$ is free as a left $\Hb$--module by 
Proposition~\ref{proof of gensit}, tensoring with $V^b$ is an exact functor. Hence,
we have the following.

\begin{Corollary} Suppose that $X$ is a right ideal of $\Hb$. Then 
$$\hatHo_{s,b}(X) = \Theta_b(X)V^b 
                = \bigoplus_{w\in\D_{(b,n-b)}}\Theta_b(X)v_bT_w, $$
where $\Theta_b(X)v_bT_w=\Theta_b(X)V^b_w\cong X$ as an $R$--module. 
In particular, if $X$ is a free $R$--module of rank $\ell$ then
$\hatHo_{s,b}(X)$ is a free $R$--module of rank $\ell\binom n b$.
\label{multiplyup}
\end{Corollary}

More generally, if 
$0\To X\To Y\To Y/X\To0$ is an exact sequence of $\Hb$--modules then
$\hatHo_{s,b}(Y/X)\cong\hatHo_{s,b}(Y)/\hatHo_{s,b}(X)
                  \cong\Theta_b(Y)V^b/\Theta_b(X)V^b$.

We are now ready to describe how the Specht modules $S^\lambda$ and
simple $\H$--modules $D^\mu$ of $\Hb$ ``multiply up'' to give
$\H$--modules.

\begin{Lemma} Suppose that $R$ is a field, $0\le b\le n$ and that
$\lambda\in\Lab$.  Then $S^\lambda\in\Modb$; moreover, if
$D^\lambda\ne(0)$ then $D^\lambda\in\Modb$.
\label{placings}

\end{Lemma}

\proof First note that by \cite[Theorem~3.7(ii)]{GL} all of the
composition factors of $S^{\lambda}$ belong to the same block;
therefore, by Theorem~\ref{Specht filtration}(i) if $\lambda\in\Lab$ then
$S^\lambda\in\Modb$. Finally, if $D^\lambda\ne(0)$ then $S^\lambda$
and $D^\lambda$ belong to the same block, so $D^\lambda$ also
belongs to $\Modb$.
\endproof

Consequently, if $\lambda\in\Lab$ then
$\hatHo_s(S^\lambda)=\hatHo_{s,b}(S^\lambda)$
and $\hatHo_s(D^\lambda)=\hatHo_{s,b}(D^\lambda)$.

Recall that $d_{\lambda\mu}=[S^\lambda:D^\mu]$ is the decomposition
multiplicity of the simple module $D^\mu$ in the Specht module
$S^\lambda$.

\begin{Proposition}\label{spechtmultplyup} Suppose that $R$ is a field and that
$f_s(q,\Q)\ne0$. Let $\lambda\in\Lab$ and $\mu\in\Lambda^+_c$ be 
multipartitions of $n$ and, as in Lemma~\ref{ulambdaactsright}, write 
$\lambda=(\sigma,\tau)$ and $\mu=(\alpha,\beta)$,  where 
$\sigma\in\Lamp(b:s)$, $\tau\in\Lamp(n-b:r-s)$,
$\alpha\in\Lamp(c:s)$ and $\beta\in\Lamp(n-c:r-s)$ (and $0\le b,c\le n$).
Then the following hold.
\begin{enumerate}\item $S^\lambda\cong\hatHo_s(S^\sigma\otimes S^\tau)
                     =\hatHo_{s,b}(S^\sigma\otimes S^\tau)$.
\item $D^\mu\cong\hatHo_s(D^\alpha\otimes D^\beta)
                =\hatHo_{s,c}(D^\alpha\otimes D^\beta)$.
Consequently, $D^\mu\ne(0)$ if and only if $D^\alpha\ne(0)$ and
$D^\beta\ne(0)$; in addition, 
$\dim D^\mu=\binom nc(\dim D^\alpha)(\dim D^\beta)$.
\item Suppose that $D^\mu\ne(0)$. Then
$$d_{\lambda\mu}=\cases d_{\sigma\alpha}d_{\tau\beta},&\If b=c,\\
                        0,&\If b\ne c.
                 \endcases$$
\end{enumerate}

\end{Proposition}

\proof For part (i) we don't actually need to assume that $R$ is a field.
By definition, $S^\sigma\otimes S^\tau=(z_\sigma\otimes z_\tau)\Hb$
where 
$z_\sigma\otimes z_\tau
    =m_\sigma\otimes m_\tau+{\bar N}^\sigma\!\!\otimes{\bar N}^\tau\!\!$.
By Corollary~\ref{multiplyup} and Lemma~\ref{nice action},
$\hatHo_{s,b}({\bar N}^\sigma\!\!\otimes{\bar N}^\tau\!\!)
   \cong\Theta_b({\bar N}^\sigma\!\!\otimes{\bar N}^\tau\!\!)V^b
   \cong V^b\Nlambar$; therefore,
$\hatHo_{s,b}(S^\sigma\otimes S^\tau)\cong S^\lambda$ by the remarks 
following Corollary~\ref{multiplyup}. (Notice that in this case the formula $\dim
S^\lambda=\binom nb\dim(S^\sigma\otimes S^\tau)$ is just the combinatorial
identity $|\Std(\lambda)|=\binom nb|\Std(\sigma)|\cdot|\Std(\tau)|$.)


For (ii)  note that because $\hatHo_{s,c}$ is an equivalence of
categories it takes simple $\Hb$--modules to simple $\H$--modules.
Hence, $D^\mu\cong\hatHo_{s,c}(D^\alpha\otimes D^\beta)$ by
induction on the dominance ordering using part (i) and
(\ref{unitriangular}); the dimension formula now follows from
Corollary~\ref{multiplyup}.  

Finally, as $\hatHo_s$ takes composition series to composition
series, part (iii) follows from (i) and (ii) and the decomposition
$\Mod_\H=\bigoplus_{a=0}^n\Mod_{\H\!(a)}$.
\endproof    

Note that part (iii) of the Proposition provides a recipe for
calculating the decomposition matrix of $\H$ from the decomposition
matrices of the ``smaller'' Ariki--Koike algebras $\H_{q,\Q_1}(b)$ and
$\H_{q,\Q_2}(n-b)$ for $0\le b\le n$.

Proposition~\ref{spechtmultplyup} and all of the other consequences of Theorem~\ref{d=2} can
be extended to the general case of Theorem~\ref{Morita equivalence}, where the
parameter set $\Q$ is partitioned into an arbitrary number of pieces. The
notation needed to describe this is rather cumbersome so we leave the
details to the reader.

\section{The cyclotomic $q$--Schur algebra}

We now extend the Morita equivalence of the previous section to
cyclotomic $q$--Schur algebras. Let $\Gamma\subseteq\Lambda$ be a
finite set of multicompositions of $n$ with the property that
whenever $\mu\in\Lambda^+$ and $\mu\gedom\lambda$ for some
$\lambda\in\Gamma$ then $\mu\in\Gamma$.  The {\sf cyclotomic
$q$--Schur algebra} associated with~$\Gamma$ is the algebra
$$\Sch_{q,\Q}(\Gamma)
        =\End_\H\Big(\bigoplus_{\lambda\in\Gamma}M^\lambda\Big).$$
For the statement of Theorem~\ref{main3} we set
$\Sch_{q,\Q}(n)=\Sch_{q,\Q}(\Lambda^+)$.

Let $\Gamma^+=\Gamma\cap\Lambda^+$ be the set of multipartitions in
$\Gamma$.  Then, by \cite[Theorem~6.6]{DJM:cyc} (see also 
(\ref{SStd basis})), the algebra $\Sch_{q,\Q}(\Gamma)$ has (cellular)
basis
$$\set{\phiST|\S\in\SStd(\lambda,\mu)\And\T\in\SStd(\lambda,\nu)
        \ForSome\mu,\nu\in\Gamma\text{and some}\lambda\in\Gamma^+},$$
where $\phiST$ is the $\H$--module homomorphism given by
$\phiST(m_\alpha h)=\delta_{\alpha\nu}m_{\S\T}h$. The basis
$\{\phiST\}$ is the {\sf semistandard basis} of
$\Sch_{q,\Q}(\Gamma)$.

As before, we fix an integer $s$, with $1\le s\le r$, such that $f_s(q,\Q)$
is invertible in $R$ and let 
$\Gamma_b=\set{\rtuple\lambda\in\Gamma|
                        b=|\lambda^{(1)}|+\dots+|\lambda^{(s)}|}.$
We need to define analogues of the sets $\Lamp(b:s)$; however, we must
be a little careful. Let
\begin{align*}
\Gaml(b:s)&=\set{\sigma\in\Lam(b:s)|(\sigma,\tau)\in\Gamma
                       \ForSome\tau\in\Lam(n-b:r-s)}
\intertext{and}
\Gamr(n-b:r-s)&=\set{\tau\in\Lam(n-b:r-s)|(\sigma,\tau)\in\Gamma
                       \ForSome\sigma\in\Lam(b:s)}.
\end{align*}
Also let $\Gamlp(b:s)=\Gaml(b:s)\cap\Lamp(b:s)$ and
$\Gamrp(n-b:r-s)=\Gamr(n-b:r-s)\cap\Lamp(n-b:r-s)$.
Consider $\Gaml(b:s)\times\Gamr(n-b:r-s)$ as a poset in the obvious way.

\begin{Lemma} Suppose that $0\le b\le n$. Then $\Gaml(b:s)\times\Gamr(n-b:r-s)$
and $\Gamma_b$ are naturally isomorphic posets. In particular, if 
$(\alpha,\beta)\in\Lamp(b:s)\times\Lamp(n-b:r-s)$ and
$(\alpha,\beta)\gedom(\sigma,\tau)$ for some
$(\sigma,\tau)\in\Gaml(b:s)\times\Gamr(n-b:r-s)$ then
$(\alpha,\beta)\in\Gaml(b:s)\times\Gamr(n-b:r-s)$.

\end{Lemma}

\proof The isomorphism is given by 
$$\((\sigma^{(1)},\dots,\sigma^{(s)}),(\tau^{(1)},\dots,\tau^{(r-s)})\)
\mapsto(\sigma^{(1)},\dots,\sigma^{(s)},\tau^{(1)},\dots,\tau^{(r-s)}).$$
This is a poset isomorphism because 
$|\sigma^{(1)}|+\dots+|\sigma^{(s)}|=b$ whenever
$(\sigma,\tau)$ is an element of $\Gaml(b:s)\times\Gamr(n-b:r-s)$.
\endproof

Let $\Sb=\Sch_{q,\Q_1}\(\Gaml(b:s)\)\otimes\Sch_{q,\Q_2}\(\Gamr(n-b:r-s)\)$,
where $\Q_1=(Q_1,\dots,Q_s)$ and $\Q_2=(Q_{s+1},\dots,Q_r)$. The point of
the Lemma is that it allows us to index the representations of $\Sb$ with
the elements of $\Gamma_b$. For convenience we identity $\Gamma_b$
and $\Gaml(b:s)\times\Gamr(n-b:r-s)$ in the sequel.

We can now give the analogue of Theorem~\ref{d=2} for a cyclotomic
$q$--Schur algebra. Rather than introduce the notation necessary to
state the general result, we consider only the special case where
$\Q$ is partitioned into two pieces. See Theorem~\ref{main3} for a special 
case of the more general result.

\begin{Theorem}\label{genmain3} Let $1\le s\le r$ and suppose that $f_s(q,\Q)$ is an
invertible element of $R$. Then the cyclotomic $q$--Schur algebra
$\Sch_{q,\Q}(\Gamma)$ is Morita equivalent to 
$$ \Sch_s(\Gamma)=\bigoplus_{b=0}^n\Sb.  $$

\end{Theorem}

Before we can give the proof we require some preparation. The basic idea is
to investigate how the functor $\Ho_s$ acts on 
$\bigoplus_{\lambda\in\Gamma}M^\lambda$; in fact, it is easier to work
with $\hatHo_s$.

\begin{Lemma} Suppose that $0\le b\le n$ and let $\lambda=(\sigma,\tau)\in\Gamma_b$.
Then $$\hatHo_{s,b}(M^\sigma\otimes M^\tau)=\theta_b(M^\lambda).$$
\label{first cut}

\end{Lemma}

\proof We first note that $u_b^+=m_{\omega_b}$ is a left factor of 
$u_\lambda^+$ because $\lambda\in\Lab$; consequently, $M^\lambda$ is a 
submodule of~$\Mb$ and so $\theta_b(M^\lambda)$ makes sense. Now
$\Theta_b(m_\sigma\otimes m_\tau)v_b=\theta_b(m_\lambda)$ by 
Lemma~\ref{nice action}; therefore, 
$\hatHo_{s,b}(M^\sigma\otimes M^\tau)=\theta_b(m_\lambda)\H
                      =\theta_b(M^\lambda)$
by Corollary~\ref{multiplyup}.
\endproof

Combining Theorem~\ref{Specht filtration} and Corollary~\ref{splitting} we see that
$\theta_b$ projects $\Mb$ onto $V^b$ and consequently that $V^b$ is
the unique direct summand of $\Mb$ which belongs to $\Modb$. Now,
$\H=\bigoplus_{c=0}^n\H(c)$ so we can write
$M^\lambda=\bigoplus_{c=0}^nM^\lambda(c)$, where $M^\lambda(c)$ is
the largest direct summand of $M^\lambda$ which is contained in
$\H(c)$.  Furthermore, since $M^\lambda$ is a quotient of $\Mb$,
$M^\lambda(c)=(0)$ whenever $c<b$ by Corollary~\ref{splitting}.  To proceed we
need to understand the direct summands $M^\lambda(c)$ of
$M^\lambda$.

Given a multicomposition $\mu$ let $\bar\mu=\rtuple{\bar\mu{}}$ be
the (unique) multipartition where $\bar\mu{}^{(c)}$ is obtained by
ordering the parts of $\mu^{(c)}$\!, for $1\le c\le r$. 

\begin{Point}{\cite[Corollary~3.5]{M:cycSchur}} There exists a family
$\set{Y^\lambda|\lambda\in\Lambda^+}$ of indecomposable
$\H$--modules which are uniquely determined by the property that,
for each $\mu\in\Lambda$,
$$M^\mu\cong Y^{\bar\mu}\oplus
            \bigoplus_{\SR{\lambda\in\Lambda^+}{\lambda\gdom\bar\mu}}
                 c_{\lambda\bar\mu}Y^\lambda$$
for some non--negative integers $c_{\lambda\bar\mu}$ $($which depend only on
$\lambda$ and $\bar\mu\,)$. Moreover,~$S^\lambda$ is a quotient of~$Y^\lambda$.
\label{indecomposables}
\end{Point}

The $Y^\lambda$ are generalizations of the {\sf Young modules} of
the symmetric groups.  For the case of the $q$--Schur algebra (that
is, when $r=1$) this result is proved in \cite{DJ:Schur}. 

Set $c_{\lambda\lambda}=1$ and $c_{\lambda\bar\mu}=0$ if
$\lambda\not\gedom\bar\mu$.  We can now identify the summand of $M^\mu$
which is contained in $\H(b)$.

\begin{Corollary} Suppose that $\mu\in\Gamma$. Then 
$\DS M^\mu(b)\cong\bigoplus_{\lambda\in\Gamma_b^+}
         c_{\lambda\bar\mu}Y^\lambda$ for $0\le b\le n$.
\label{second cut}

\end{Corollary}

\proof If $\lambda\in\Lambda_b^+$ then $S^\lambda\in\Modb$ by
Lemma~\ref{placings}; therefore, $Y^\lambda$ also belongs to $\Modb$, since
$Y^\lambda$ is indecomposable and $S^\lambda$ is a quotient of $Y^\lambda$.
(Observe that $\bar\mu\gedom\mu$; consequently, if $\lambda\in\Lambda^+$
and $\lambda\gedom\bar\mu$ then $\lambda\gedom\mu$ and so
$\lambda\in\Gamma^+$.) As the decomposition of $M^\mu$ into a direct sum of
indecomposables is unique up to isomorphism, the result follows. 
\endproof

Let $M_\Gamma=\bigoplus_{\mu\in\Gamma}M^\mu$ and
$M_{\Gamma_b}=\bigoplus_{(\sigma,\tau)\in\Gamma_b}
                  M^\sigma\otimes M^\tau$.
Then, by definition, $\Sch_{q,\Q}(\Gamma)=\End_\H(M_\Gamma)$ and
$\Sch_s(\Gamma)=\bigoplus_{b=0}^n\End_{\Hb}(M_{\Gamma_b})$. Finally,
we set $M_{\Gamma,s}=\bigoplus_{b=0}^n\hatHo_{s,b}(M_{\Gamma_b})$.
These modules are the analogues of {\sf $q$--tensor space} for the
various cyclotomic $q$--Schur algebras; compare with
\cite{DJ:qWeyl}.

\begin{Lemma} Suppose that $f_s(q,\Q)$ is invertible. Then
$M_{\Gamma,s}$ and $M_\Gamma$ have the same
set of indecomposable direct summands.

\end{Lemma}

\proof Applying the definitions,
\begin{align*}
M_{\Gamma,s} &=\bigoplus_{b=0}^n\hatHo_{\s,b}(M_{\Gamma_b})
     \cong\bigoplus_{b=0}^n\bigoplus_{(\sigma,\tau)\in\Gamma_b}
          \hatHo_{\s,b}(M^\sigma\otimes M^\tau).
\intertext{Now, if $\mu=(\sigma,\tau)\in\Gamma_b$ then
$\hatHo_{\s,b}(M^\sigma\otimes M^\tau)=M^\mu(b)$ by Lemma~\ref{first cut};
on the other hand,
$M^\mu(b)\cong\bigoplus_{\lambda\in\Gamma_b^+}c_{\lambda\bar\mu}
                   Y^\lambda$
by Corollary~\ref{second cut}. Therefore,}
M_{\Gamma,s} &\cong\bigoplus_{b=0}^n\bigoplus_{\mu\in\Gamma_b}
        \bigoplus_{\lambda\in\Gamma_b^+}c_{\lambda\bar\mu}Y^\lambda.
\end{align*}
Consequently, $\set{Y^\lambda|\lambda\in\Gamma^+}$ is a complete set of
isomorphism classes of indecomposable direct summands of $M_{\Gamma,s}$.
By~(\ref{indecomposables}) this is also the set of indecomposable direct
summands of $M_\Gamma$, so this proves the Lemma.
\endproof

\proofof{genmain3} By \Last, the sets of indecomposable direct
summands of the $\H$--modules $M_\Gamma$ and $M_{\Gamma,s}$ coincide,
except that the multiplicities of each $Y^\lambda$ in the two modules will
typically differ. Therefore, by a general argument, the $\H$--module
endomorphism rings of $M_\Gamma$ and $M_{\Gamma,s}$ are Morita equivalent.
On the other hand, since~$\hatHo_s$ is an equivalence of categories, it
induces an isomorphism of endomorphism rings;~so,
$$\Sch_s(\Gamma)=\bigoplus_{b=0}^n\End_{\Hb}(M_{\Gamma_b})
            \cong\bigoplus_{b=0}^n\End_\H\(\hatHo_{s,b}(M_{\Gamma_b})\)
            \cong\End_\H\(M_{\Gamma,s}\),$$ 
where the last isomorphism is a consequence of the decomposition
$\Mod_\H=\bigoplus_{b=0}^n\Modb$. As we have already observed that
$\End_\H\(M_{\Gamma,s}\)$ is Morita equivalent
to~$\Sch_{q,\Q}(\Gamma)=\End_\H(M_\Gamma)$, the proof of Theorem~\ref{genmain3} 
is complete.
\endproof

To conclude, we note that Proposition~\ref{spechtmultplyup} also generalizes to
the cyclotomic $q$--Schur algebra case. Formally, the proof is similar to 
the Ariki--Koike case, so we omit the details. 

Recall that the cell modules of $\Sch_{q,\Q}(\Gamma)$ are called
{\sf Weyl modules} and are denoted by $W^\lambda$, for
$\lambda\in\Gamma^+$. As with the Specht modules, there is a
symmetric associative bilinear form on $W^\lambda$; let
$F^\lambda=W^\lambda/\rad W^\lambda$, where $\rad W^\lambda$ is the
radical of this form. Then $\set{F^\lambda|\lambda\in\Gamma^+}$ is a
complete set of pairwise non--isomorphic irreducible
$\Sch_{q,\Q}(\Gamma)$--modules by \cite[Theorem~6.16]{DJM:cyc}.
Furthermore, by \cite[Theorem~2.3]{M:cycSchur} if $D^\mu\ne(0)$ then
$[W^\lambda:F^\mu]=[S^\lambda:D^\mu]$; because of this we abuse
notation and write $d_{\lambda\mu}=[W^\lambda:F^\mu]$ for
all $\lambda,\mu\in\Gamma^+$.

Again, we write $\hatHo_s$ and $\hatHo_{s,b}$ for the induced functors
between the categories $\Mod_{\Sch_s(\Gamma)}$ and 
$\Mod_{\Sch_{q,\Q}(\Gamma)}$ and their subcategories.

\begin{Corollary} Suppose that $R$ is a field and that $f_s(q,\Q)\ne0$. Suppose
that $\lambda=(\sigma,\tau)\in\Gamma_b^+$ and
$\mu=(\alpha,\beta)\in\Gamma^+_c$ where $0\le b,c\le n$.  Then we 
have the following.
\begin{enumerate}\item $W^\lambda\cong\hatHo_s(W^\sigma\otimes W^\tau)
                     =\hatHo_{s,b}(W^\sigma\otimes W^\tau)$.
\item $F^\mu\cong\hatHo_s(F^\alpha\otimes F^\beta)
                =\hatHo_{s,c}(F^\alpha\otimes F^\beta)$.
Consequently, $F^\mu\ne(0)$ if and only if $F^\alpha\ne(0)$ and
$F^\beta\ne(0)$.
\item Finally, if $b=c$ then 
$d_{\lambda\mu}= d_{\sigma\alpha}d_{\tau\beta}$; otherwise,
$d_{\lambda\mu}=0$.
\end{enumerate}
\end{Corollary}

\let\em\it

\makeatletter\@address
\end{document}